\documentclass{amsproc}
\begin{document}
\title[Geometric Crystals on Flag Varieties and 
Unipotent Radicals]{Geometric Crystals on 
Flag Varieties and  Unipotent Subgroups of Classical Groups}

\author{Mana I\textsc{garashi}}
\address{M.I.:\,Department of Mathematics, Sophia University, 
Kioicho 7-1, Chiyoda-ku, Tokyo 102-8554, Japan}
\email{mana-i@hoffman.cc.sophia.ac.jp}
\author{Toshiki N\textsc{akashima}}
\address{T.N.:\,Department of Mathematics, Sophia University, 
Kioicho 7-1, Chiyoda-ku, Tokyo 102-8554, Japan}
\email{toshiki@mm.sophia.ac.jp,\,\,toshiki@sophia.ac.jp}
\thanks{The author was supported in part by  JSPS Grants in 
Aid for Scientific Research \#19540050.}

\date{January , 2005}


\begin{abstract}
For a classical simple algebraic group $G$ we obtain the affirmative answer  
for the conjecture in \cite{N2} 
that 
there exists an isomorphism between the geometric crystal
on the flag variety and the one on the unipotent subgroup $U^-$.
\end{abstract}

\maketitle

\renewcommand{\labelenumi}{$($\roman{enumi}$)$}
\renewcommand{\labelenumii}{$(${\rm \alph{enumii}}$)$}
\font\germ=eufm10

\newcommand{\cI}{{\mathcal I}}
\newcommand{\cA}{{\mathcal A}}
\newcommand{\cB}{{\mathcal B}}
\newcommand{\cC}{{\mathcal C}}
\newcommand{\cD}{{\mathcal D}}
\newcommand{\cF}{{\mathcal F}}
\newcommand{\cH}{{\mathcal H}}
\newcommand{\cK}{{\mathcal K}}
\newcommand{\cL}{{\mathcal L}}
\newcommand{\cM}{{\mathcal M}}
\newcommand{\cN}{{\mathcal N}}
\newcommand{\cO}{{\mathcal O}}
\newcommand{\cS}{{\mathcal S}}
\newcommand{\cT}{{\mathcal T}}
\newcommand{\cV}{{\mathcal V}}
\newcommand{\fra}{\mathfrak a}
\newcommand{\frb}{\mathfrak b}
\newcommand{\frc}{\mathfrak c}
\newcommand{\frd}{\mathfrak d}
\newcommand{\fre}{\mathfrak e}
\newcommand{\frf}{\mathfrak f}
\newcommand{\frg}{\mathfrak g}
\newcommand{\frh}{\mathfrak h}
\newcommand{\fri}{\mathfrak i}
\newcommand{\frj}{\mathfrak j}
\newcommand{\frk}{\mathfrak k}
\newcommand{\frI}{\mathfrak I}
\newcommand{\fm}{\mathfrak m}
\newcommand{\frn}{\mathfrak n}
\newcommand{\frp}{\mathfrak p}
\newcommand{\fq}{\mathfrak q}
\newcommand{\frr}{\mathfrak r}
\newcommand{\frs}{\mathfrak s}
\newcommand{\frt}{\mathfrak t}
\newcommand{\fru}{\mathfrak u}
\newcommand{\frA}{\mathfrak A}
\newcommand{\frB}{\mathfrak B}
\newcommand{\frF}{\mathfrak F}
\newcommand{\frG}{\mathfrak G}
\newcommand{\frH}{\mathfrak H}
\newcommand{\frJ}{\mathfrak J}
\newcommand{\frN}{\mathfrak N}
\newcommand{\frP}{\mathfrak P}
\newcommand{\frT}{\mathfrak T}
\newcommand{\frU}{\mathfrak U}
\newcommand{\frV}{\mathfrak V}
\newcommand{\frX}{\mathfrak X}
\newcommand{\frY}{\mathfrak Y}
\newcommand{\frZ}{\mathfrak Z}
\newcommand{\rA}{\mathrm{A}}
\newcommand{\rC}{\mathrm{C}}
\newcommand{\rd}{\mathrm{d}}
\newcommand{\rB}{\mathrm{B}}
\newcommand{\rD}{\mathrm{D}}
\newcommand{\rE}{\mathrm{E}}
\newcommand{\rH}{\mathrm{H}}
\newcommand{\rK}{\mathrm{K}}
\newcommand{\rL}{\mathrm{L}}
\newcommand{\rM}{\mathrm{M}}
\newcommand{\rN}{\mathrm{N}}
\newcommand{\rR}{\mathrm{R}}
\newcommand{\rT}{\mathrm{T}}
\newcommand{\rZ}{\mathrm{Z}}
\newcommand{\bbA}{\mathbb A}
\newcommand{\bbC}{\mathbb C}
\newcommand{\bbG}{\mathbb G}
\newcommand{\bbF}{\mathbb F}
\newcommand{\bbH}{\mathbb H}
\newcommand{\bbP}{\mathbb P}
\newcommand{\bbN}{\mathbb N}
\newcommand{\bbQ}{\mathbb Q}
\newcommand{\bbR}{\mathbb R}
\newcommand{\bbV}{\mathbb V}
\newcommand{\bbZ}{\mathbb Z}
\newcommand{\adj}{\operatorname{adj}}
\newcommand{\Ad}{\mathrm{Ad}}
\newcommand{\Ann}{\mathrm{Ann}}
\newcommand{\rcris}{\mathrm{cris}}
\newcommand{\ch}{\mathrm{ch}}
\newcommand{\cl}{\colon}
\newcommand{\coker}{\mathrm{coker}}
\newcommand{\diag}{\mathrm{diag}}
\newcommand{\Diff}{\mathrm{Diff}}
\newcommand{\Dist}{\mathrm{Dist}}
\newcommand{\rDR}{\mathrm{DR}}
\newcommand{\ev}{\mathrm{ev}}
\newcommand{\Ext}{\mathrm{Ext}}
\newcommand{\cExt}{\mathcal{E}xt}
\newcommand{\fin}{\mathrm{fin}}
\newcommand{\Frac}{\mathrm{Frac}}
\newcommand{\GL}{\mathrm{GL}}
\newcommand{\Hom}{\mathrm{Hom}}
\newcommand{\hd}{\mathrm{hd}}
\newcommand{\rht}{\mathrm{ht}}
\newcommand{\id}{\mathrm{id}}
\newcommand{\im}{\mathrm{im}}
\newcommand{\inc}{\mathrm{inc}}
\newcommand{\ind}{\mathrm{ind}}
\newcommand{\coind}{\mathrm{coind}}
\newcommand{\Lie}{\mathrm{Lie}}
\newcommand{\Max}{\mathrm{Max}}
\newcommand{\mult}{\mathrm{mult}}
\newcommand{\op}{\mathrm{op}}
\newcommand{\ord}{\mathrm{ord}}
\newcommand{\pt}{\mathrm{pt}}
\newcommand{\qt}{\mathrm{qt}}
\newcommand{\rad}{\mathrm{rad}}
\newcommand{\res}{\mathrm{res}}
\newcommand{\rgt}{\mathrm{rgt}}
\newcommand{\rk}{\mathrm{rk}}
\newcommand{\SL}{\mathrm{SL}}
\newcommand{\seteq}{\mathbin{:=}}
\newcommand{\soc}{\mathrm{soc}}
\newcommand{\Spec}{\mathrm{Spec}}
\newcommand{\St}{\mathrm{St}}
\newcommand{\supp}{\mathrm{supp}}
\newcommand{\Tor}{\mathrm{Tor}}
\newcommand{\Tr}{\mathrm{Tr}}
\newcommand{\wt}{\mathrm{wt}}
\newcommand{\Ab}{\mathbf{Ab}}
\newcommand{\Alg}{\mathbf{Alg}}
\newcommand{\Grp}{\mathbf{Grp}}
\newcommand{\Mod}{\mathbf{Mod}}
\newcommand{\Sch}{\mathbf{Sch}}\newcommand{\bfmod}{{\bf mod}}
\newcommand{\Qc}{\mathbf{Qc}}
\newcommand{\Rng}{\mathbf{Rng}}
\newcommand{\Top}{\mathbf{Top}}
\newcommand{\Var}{\mathbf{Var}}
\newcommand{\gromega}{\langle\omega\rangle}
\newcommand{\lbr}{\begin{bmatrix}}
\newcommand{\rbr}{\end{bmatrix}}
\newcommand{\forb}{\bigcirc\kern-2.8ex \because}
\newcommand{\forbb}{\bigcirc\kern-3.0ex \because}
\newcommand{\forbbb}{\bigcirc\kern-3.1ex \because}
\newcommand{\cd}{commutative diagram }
\newcommand{\SpS}{spectral sequence}
\newcommand\C{\mathbb C}
\newcommand\hh{{\hat{H}}}
\newcommand\eh{{\hat{E}}}
\newcommand\F{\mathbb F}
\newcommand\fh{{\hat{F}}}

\def\AA{{\mathcal A}}
\def\al{\alpha}
\def\bq{B_q(\ge)}
\def\bqm{B_q^-(\ge)}
\def\bqz{B_q^0(\ge)}
\def\bqp{B_q^+(\ge)}
\def\beneme{\begin{enumerate}}
\def\beq{\begin{equation}}
\def\beqn{\begin{eqnarray}}
\def\beqnn{\begin{eqnarray*}}
\def\bigsl{{\hbox{\fontD \char'54}}}
\def\bbra#1,#2,#3{\left\{\begin{array}{c}\hspace{-5pt}
#1;#2\\ \hspace{-5pt}#3\end{array}\hspace{-5pt}\right\}}
\def\cd{\cdots}
\def\CC{\hbox{\bf C}}
\def\ddd{\hbox{\germ D}}
\def\del{\delta}
\def\Del{\Delta}
\def\Delr{\Delta^{(r)}}
\def\Dell{\Delta^{(l)}}
\def\Delb{\Delta^{(b)}}
\def\Deli{\Delta^{(i)}}
\def\Delre{\Delta^{\rm re}}
\def\ei{e_i}
\def\eit{\tilde{e}_i}
\def\eneme{\end{enumerate}}
\def\ep{\epsilon}
\def\eeq{\end{equation}}
\def\eeqn{\end{eqnarray}}
\def\eeqnn{\end{eqnarray*}}
\def\fit{\tilde{f}_i}
\def\FF{{\rm F}}
\def\ft{\tilde{f}}
\def\gau#1,#2{\left[\begin{array}{c}\hspace{-5pt}#1\\
\hspace{-5pt}#2\end{array}\hspace{-5pt}\right]}
\def\ge{\hbox{\germ g}}
\def\gl{\hbox{\germ gl}}
\def\hom{{\hbox{Hom}}}
\def\ify{\infty}
\def\io{\iota}
\def\kp{k^{(+)}}
\def\km{k^{(-)}}
\def\llra{\relbar\joinrel\relbar\joinrel\relbar\joinrel\rightarrow}
\def\lan{\langle}
\def\lar{\longrightarrow}
\def\max{{\rm max}}
\def\lm{\lambda}
\def\Lm{\Lambda}
\def\mapright#1{\smash{\mathop{\longrightarrow}\limits^{#1}}}
\def\hookright#1{\smash{\mathop{\hookrightarrow}\limits^{#1}}}
\def\mm{{\bf{\rm m}}}
\def\nd{\noindent}
\def\nn{\nonumber}
\def\nnn{\hbox{\germ n}}
\def\catob{{\mathcal O}(B)}
\def\oint{{\mathcal O}_{\rm int}(\ge)}
\def\ot{\otimes}
\def\op{\oplus}
\def\opi{\ovl\pi_{\lm}}
\def\ovl{\overline}
\def\plm{\Psi^{(\lm)}_{\io}}
\def\qq{\qquad}
\def\q{\quad}
\def\qed{\hfill\framebox[2mm]{}}
\def\QQ{\hbox{\bf Q}}
\def\qi{q_i}
\def\qii{q_i^{-1}}
\def\ra{\rightarrow}
\def\ran{\rangle}
\def\rlm{r_{\lm}}
\def\ssl{\mathfrak{sl}}
\def\slh{\widehat{\ssl_2}}
\def\ge{\hbox{\germ g}}
\def\ti{t_i}
\def\tii{t_i^{-1}}
\def\til{\tilde}
\def\tm{\times}
\def\tt{{\hbox{\germ{t}}}}
\def\ttt{\hbox{\germ t}}
\def\ua{U_{\AA}}
\def\ue{U_{\vep}}
\def\uq{U_q(\ge)}
\def\ufin{U^{\rm fin}_{\vep}}
\def\ufinp{(U^{\rm fin}_{\vep})^+}
\def\ufinm{(U^{\rm fin}_{\vep})^-}
\def\ufinz{(U^{\rm fin}_{\vep})^0}
\def\uqm{U^-_q(\ge)}
\def\uqp{U^+_q(\ge)}
\def\uqmq{{U^-_q(\ge)}_{\bf Q}}
\def\uqpm{U^{\pm}_q(\ge)}
\def\uqq{U_{\bf Q}^-(\ge)}
\def\uqz{U^-_{\bf Z}(\ge)}
\def\ures{U^{\rm res}_{\AA}}
\def\urese{U^{\rm res}_{\vep}}
\def\uresez{U^{\rm res}_{\vep,\ZZ}}
\def\util{\widetilde\uq}
\def\uup{U^{\geq}}
\def\ulow{U^{\leq}}
\def\bup{B^{\geq}}
\def\blow{\ovl B^{\leq}}
\def\vep{\varepsilon}
\def\vp{\varphi}
\def\vpi{\varphi^{-1}}
\def\VV{{\mathcal V}}
\def\xii{\xi^{(i)}}
\def\Xiioi{\Xi_{\io}^{(i)}}
\def\WW{{\mathcal W}}
\def\wtil{\widetilde}
\def\what{\widehat}
\def\wpi{\widehat\pi_{\lm}}
\def\ZZ{\mathbb Z}
\def\spsp(#1,#2){\begin{pmatrix}
#1,\\ \hline#2\end{pmatrix}}
\def\sps(#1;#2){(#1 | #2)}

\theoremstyle{definition}
\newtheorem{df}{\bf Definition}[section]
\theoremstyle{plain}
\newtheorem{pro}[df]{\bf Proposition}
\newtheorem{lem}[df]{\bf Lemma}
\newtheorem{thm}[df]{\bf Theorem}
\newtheorem{cor}[df]{\bf Corollary}

\renewcommand{\thesection}{\arabic{section}}
\section{Introduction}
\setcounter{equation}{0}
\renewcommand{\theequation}{\thesection.\arabic{equation}}

The theory of geometric crystal for semi-simple case 
has been introduced in \cite{BK}
as an geometric analogue of Kashiwara's crystal theory.
In \cite{N} it has been extended to Kac-Moody setting and 
the geometric
crystals on Schubert variety $\ovl X_w$ has been 
introduced therein, where $w$ is a Weyl group
element. In \cite{N2} we constructed geometric crystals on 
the unipotent radical $U^-\subset B^-$ of 
a semi-simple algebraic group $G$, where $B^-$ is an opposite Borel 
subgroup and showed that in the case $G=SL_n(\bbC)$ 
it is isomorphic to the geometric crystal on the flag variety $X=\ovl X_{w_0}$
where $w_0$ is the longest element in the corresponding Weyl group.
In \cite{N2} we conjectured  that 
for any semi-simple case there exists such 
isomorphism and in this article we obtained the isomorphism
between $U^-$ and the flag variety 
for classical simple algebraic groups.

Here we explain more details.
Let $B^\pm\subset G$ be the Borel subgroups and $U^\pm$ 
their unipotent radicals.
As mentioned above, in \cite{N} we constructed geometric crystals on Schubert
varieties, whose dimension is finite. 
Nevertheless, we can not apply the method in 
\cite{N} to the full flag variety
since it is infinite dimensional for general Kac-Moody cases.
Thus, we considered alternative way to obtain geometric crystal
structure on the opposite unipotent radical $U^-\subset B^-$ 
which is birationally isomorphic to the full flag variety 
$X=G/B$.

A variety $V$ is called a unipotent crystal if it has a rational 
$U$-action 
and there exists a rational map from $V$ to the opposite 
Borel subgroup $B^-$ commuting with the $U$-action, 
where $B^-$ is equipped with the canonical rational $U$-action by:
\[
U\times B^-\,\,\hookright{m}\,\, G\mapright{\sim}B^-\times
U\mapright{proj}B^-.
\]
One of the most crucial properties of unipotent crystals is that 
certain geometric crystal is induced from a 
unipotent crystal canonically (see \ref{u->g}). 

We introduce a criterion for the existence of unipotent crystal 
on $U^-$ in \cite[Lemma 3.2 ]{N2} (see also Lemma \ref{suf} below.), 
which is applicable to general Kac-Moody cases though 
it is applied to only simple cases in this article.
Let us explain the criterion more precisely.
To obtain the unipotent crystal structure on $U^-$, it is required to
get certain rational map ${\mathcal T}:U^-\to T$ with the properties:
for $x\in U$ and $u\in U^-$
\[
 {\mathcal T}(\pi^{--}(xu))=\pi^0(xu) {\mathcal T}(u).
\]
Then defining a morphism $\cF:U^-\to B^-$ by $\cF(u):=u {\mathcal T}(u)$, 
$\cF$ becomes a $U$-morphism and we obtain the unipotent crystal
structure on $U^-$. 
To realize the above $\cT$ in this article we construct rational 
functions $\{F^{(n)}_i\}_{i=1,\cd,n}$ on $U^-$, 
each of which is defined as a matrix
element in the fundamental representation $L(\Lm_i)$ 
and possesses some special properties, where $\Lm_i$ is
the $i$-th fundamental weight of $\ge={\rm Lie}(G)$.
Using this rational functions $F^{(n)}_i$, we define
$\cT(u)=\prod_i\al_i^\vee(F_i^{(n)}(u)^{-1})$, which satisfies the
criterion (see Sect.4.) and then we have a unipotent crystal and the
induced geometric crystal on $U^-$.

The crucial task in this paper is writing down the explicit forms 
of the function $F^{(n)}_i$ $(i=1,\cd,n)$.
Then, using them we can check that 
there exists an isomorphism between geometric crystal on the flag
variety $X$ and the one on $U^{-}$ for the types $A_n,B_n,C_n$
and $D_n$.

Since in \cite{N2} we have made several typographical errors 
in the proof of Lemma 3.2
and modified the definition of the function $F^{(n)}_i$, we shall give
the proof of Lemma \ref{suf} in this article 
and introduce the bilinear form on 
irreducible highest weight $\ge$ module in order to 
redefine $F^{(n)}_i$.

In the last section, we give a conjecture that for all semi-simple cases
there would exist
an isomorphism of geometric crystals between a Schubert 
variety $\ovl X_w$ and $U^-_\io$ which is a dense subset in $U^-$ 
associated with a reduced word $\io$ of $w$. 
Indeed, the result in this paper would be a part of this conjecture 
for the case $w=w_0$ the longest element. 

Though our method constructing unipotent crystal structure 
on $U^-$ is valid for arbitrary Kac-Moody setting, in this
article we only treated simple cases.
Our further aim is to apply this to affine Kac-Moody cases, more
details, to find certain good functions like $F^{(n)}_i$'s for affine
cases, which would be expected to be interesting and important
from the view point of representation theory of affine Kac-Moody algebras.

\renewcommand{\thesection}{\arabic{section}}
\section{Geometric Crystals and Unipotent Crystals}\label{2}
\setcounter{equation}{0}
\renewcommand{\theequation}{\thesection.\arabic{equation}}


The notations and definitions here follow
\cite{Kac,KNO, KP, Ku2, N, N2}.
\subsection{Geometric Crystals}
\label{KM}
Fix a symmetrizable generalized Cartan matrix
 $A=({\bf a}_{ij})_{i,j\in I}$, where $I$ is a finite index set.
Let $(\tt,\{\al_i\}_{i\in I},\{h_i\}_{i\in I})$ 
be the associated
root data, where ${\tt}$ is the vector space 
over $\bbC$ with 
dimension $|I|+$ corank$(A)$, and 
$\{\al_i\}_{i\in I}\subset\tt^*$ and 
$\{h_i\}_{i\in I}\subset\tt$
are linearly independent indexed sets 
satisfying $\al_j(h_i)={\bf a}_{ij}$.

The Kac-Moody Lie algebra $\ge=\ge(A)$ associated with $A$
is the Lie algebra over $\bbC$ generated by $\tt$, the 
Chevalley generators $e_i$ and $f_i$ $(i\in I)$
with the usual defining relations (\cite{KP},\cite{PK}).
Note that if $A$ is a Cartan matrix, the corresponding Lie 
algebra $\ge$ is a semi-simple complex Lie algebra.
There is the root space decomposition 
$\ge=\bigoplus_{\al\in \tt^*}\ge_{\al}$.
Denote the set of roots by 
$\Delta:=\{\al\in \tt^*|\al\ne0,\,\,\ge_{\al}\ne(0)\}$.
Set $Q=\sum_i\bbZ \al_i$, $Q_+=\sum_i\bbZ_{\geq0} \al_i$
and $\Delta_+:=\Delta\cap Q_+$.
An element of $\Delta_+$ is called a positive root.
Define simple reflections $s_i\in{\rm Aut}(\tt)$ $(i\in I)$ by
$s_i(h):=h-\al_i(h)h_i$, which generate the Weyl group $W$.
We also define the action of $W$ on $\tt^*$ by
$s_i(\lm):=\lm-\lm(h_i)\al_i$.
Set $\Delre:=\{w(\al_i)|w\in W,\,\,i\in I\}$.

Let $\ge'$ be the derived Lie algebra of $\ge$ and 
$G$ the Kac-Moody group associated 
with $\ge'$(\cite{PK}).
Let $U_{\al}:=\exp\ge_{\al}$ $(\al\in \Delre)$
be an one-parameter subgroup of $G$.
The group $G$ is generated by $U_{\al}$ $(\al\in \Delre)$.
Let $U^{\pm}$ be the subgroups generated by $U_{\pm\al}$
($\al\in \Delre_+=\Delre\cap Q_+$), {\it i.e.,}
$U^{\pm}:=\lan U_{\pm\al}|\al\in\Del^{\rm re}_+\ran$.

For any $i\in I$, there exists a unique homomorphism;
$\phi_i:SL_2(\bbC)\rightarrow G$ such that
\[
\phi_i\left(
\left(
\begin{array}{cc}
1&t\\
0&1
\end{array}
\right)\right)=\exp t e_i,\,\,
 \phi_i\left(
\left(
\begin{array}{cc}
1&0\\
t&1
\end{array}
\right)\right)=\exp t f_i\,(t\in\bbC).
\]
Set $x_i(t):=\exp{t e_i}$, $y_i(t):=\exp{t f_i}$, 
$T_i:=\phi_i(\{{\rm diag}(t,t^{-1})|t\in\bbC\})$ and 
$N_i:=N_{G_i}(T_i)$. Let
$T$ (resp. $N$) be the subgroup of $G$ generated by $T_i$
(resp. $N_i$), which is called a {\it maximal torus} in $G$ and
$B^{\pm}=U^{\pm}T$ be the Borel subgroup of $G$.
We have the isomorphism
$\phi:W\mapright{\sim}N/T$ defined by $\phi(s_i)=N_iT/T$.
An element $\ovl s_i:=x_i(-1)y_i(1)x_i(-1)$ is in 
$N_G(T)$, which is a representative of 
$s_i\in W=N_G(T)/T$.

\begin{df}
\label{def-gc}
Let $X$ be an ind-variety over $\bbC$, $\gamma_i$ and $\vep_i$ 
$(i\in I)$ rational functions on $X$, and 
$e_i:\bbC^\times\times X\to X$ a rational $\bbC^\times$-action.
A quadruple $(X,\{e_i\}_{i\in I},\{\gamma_i,\}_{i\in I},
\{\vep_i\}_{i\in I})$ is a 
$G$ (or $\ge$)-{\it geometric crystal} 
if
\begin{enumerate}
\item
$(\{1\}\times X)\cap dom(e_i)$ 
is open dense in $\{1\}\times X$ for any $i\in I$, where
$dom(e_i)$ is the domain of definition of
$e_i\cl\C^\times\times X\to X$.
\item
The rational functions  $\{\gamma_i\}_{i\in I}$ satisfy
$\gamma_j(e^c_i(x))=c^{{\bf a}_{ij}}\gamma_j(x)$ for any $i,j\in I$.
\item
$e_i$ and $e_j$ satisfy the following relations:
\[
 \begin{array}{lll}
&e^{c_1}_{i}e^{c_2}_{j}
=e^{c_2}_{j}e^{c_1}_{i}&
{\rm if }\,{\bf a}_{ij}={\bf a}_{ji}=0,\label{0}\\
&e^{c_1}_{i}e^{c_1c_2}_{j}e^{c_2}_{i}
=e^{c_2}_{j}e^{c_1c_2}_{i}e^{c_1}_{j}&
{\rm if }\,{\bf a}_{ij}={\bf a}_{ji}=-1,
\\
&
e^{c_1}_{i}e^{c^2_1c_2}_{j}e^{c_1c_2}_{i}e^{c_2}_{j}
=e^{c_2}_{j}e^{c_1c_2}_{i}e^{c^2_1c_2}_{j}e^{c_1}_{i}&
{\rm if }\,{\bf a}_{ij}=-2,\,
{\bf a}_{ji}=-1,
\\
&
e^{c_1}_{i}e^{c^3_1c_2}_{j}e^{c^2_1c_2}_{i}
e^{c^3_1c^2_2}_{j}e^{c_1c_2}_{i}e^{c_2}_{j}
=e^{c_2}_{j}e^{c_1c_2}_{i}e^{c^3_1c^2_2}_{j}e^{c^2_1c_2}_{i}
e^{c^3_1c_2}_je^{c_1}_i&
{\rm if }\,{\bf a}_{ij}=-3,\,{\bf a}_{ji}=-1,
\end{array}
\]
\item
The rational functions $\{\vep_i\}_{i\in I}$ satisfy
$\vep_i(e_i^c(x))=c^{-1}\vep_i(x)$ and 
$\vep_i(e_j^c(x))=\vep_i(x)$ if $a_{i,j}=a_{j,i}=0$.
\end{enumerate}
\end{df}
The relations in (iii) is called 
{\it Verma relations}.
If $\chi=(X,\{e_i\},\,\{\gamma_i\},\{\vep_i\})$
satisfies the conditions (i), (ii) and (iv), 
we call $\chi$ a {\it pre-geometric crystal}.

\nd
{\bf Remark.}
The last condition (iv) is slightly modified from 
\cite{KNO, N, N2, N3, N5} since all $\vep_i$ appearing in these references 
satisfy the new condition and this condition is required 
to define ''epsilon systems'' (\cite{N6}).

\subsection{Unipotent Crystals}\label{uni-cry}

In the sequel, we denote the unipotent subgroup 
$U^+$ by $U$. 
We define unipotent crystals (see \cite{BK},\cite{N}) 
associated to Kac-Moody groups. 
\begin{df}
Let $X$ be an ind-variety over $\bbC$ and 
$\al:U\times X\rightarrow X$ a rational $U$-action
such that $\al$ is defined on $\{e\}\times X$. Then, 
the pair ${\bf X}=(X,\al)$ is called a $U$-{\it variety}. 
For $U$-varieties ${\bf X}=(X,\al_X)$
and ${\bf Y}=(Y,\al_Y)$, 
a rational map
$f:X\rightarrow Y$ is called a 
$U$-{\it morphism} if it commutes with
the action of $U$.
\end{df}
Now, we define a $U$-variety structure on $B^-=U^-T$.
As in \cite{Ku2},  the Borel subgroup $B^-$ is an ind-subgroup of $G$ and hence
an ind-variety over $\bbC$.
The multiplication map in $G$ induces the open embedding;
$ B^-\times U\hookrightarrow G,$
which is a birational map. 
Let us denote the inverse birational map by 
$ g:G\longrightarrow B^-\times U$.
Then we define the rational maps 
$\pi^-:G\rightarrow B^-$ and 
$\pi:G\rightarrow U$ by 
$\pi^-:={\rm proj}_{B^-}\circ g$ 
and $\pi:={\rm proj}_U\circ g$.
Now we define the rational $U$-action $\al_{B^-}$ on $B^-$ by 
\[
 \al_{B^-}:=\pi^-\circ m:U\times B^-\longrightarrow B^-,
\]
where $m$ is the multiplication map in $G$.
Then we get $U$-variety ${\bf B}^-=(B^-,\al_{B^-})$.
\begin{df}
\label{uni-def}
\begin{enumerate}
\item
Let ${\bf X}=(X,\al)$ 
be a $U$-variety and $f:X \rightarrow {\bf B^-}$ 
a $U$-morphism.
The pair $({\bf X}, f)$ is called 
a {\it unipotent $G$-crystal}
or, for short, {\it unipotent crystal}.
\item
Let $({\bf X},f_X)$ and $({\bf Y},f_Y)$ 
be unipotent crystals.
A $U$-morphism $g:\bf X\ra \bf Y$ is called a {\it morphism of 
unipotent crystals} if $f_X=f_Y\circ g$.
In particular, if $g$ is a birational map
of ind-varieties, it is called an {\it isomorphism of 
unipotent crystals}.
\end{enumerate}
\end{df}
We define a product of 
unipotent crystals following \cite{BK}.
For unipotent crystals $({\bf X},f_X)$, $({\bf Y},f_Y)$, 
define a morphism 
$\al_{X\times Y}:U\tm X\tm Y\rightarrow X\tm Y$ by
\begin{equation}
\al_{X\tm Y}(u,x,y):=(\al_X(u,x),\al_Y(\pi(u\cdot f_X(x)),y)).
\label{XY}
\end{equation}
\begin{thm}[\cite{BK}]
\label{prod}
\hspace{-30pt}
\begin{enumerate}
\item
The morphism $\al_{X\tm Y}$ defined above 
is a rational $U$-action
on $X\tm Y$.
\item
Let ${\bf m}:B^-\tm B^-\rightarrow B^-$ 
be a multiplication map 
and $f=f_{X\tm Y}:X\tm Y\rightarrow B^-$ be the 
rational map defined by 
\[
f_{X\tm Y}:={\bf m}\circ( f_X\tm f_Y).
\]
Then $ f_{X\tm Y}$ is a $U$-morphism and 
$({\bf X\tm Y}, f_{X\tm Y})$ is a unipotent crystal, 
which we call a product of unipotent crystals  
$({\bf X},f_X)$ and $({\bf Y},f_Y)$.
\item 
Product of unipotent crystals is associative.
\end{enumerate}
\end{thm}

\subsection{From Unipotent Crystals to Geometric Crystals}
\label{u->g}
For $i\in I$, 
set $U^\pm_i:=U^\pm\cap \bar s_i U^\mp\bar s_i^{-1}$ and
$U_\pm^i:=U^\pm\cap \bar s_i U^\pm\bar s_i^{-1}$.
Indeed, $U^\pm_i=U_{\pm \al_i}$.
Set 
\[
 Y_{\pm\al_i}:=
\lan x_{\pm\al_i}(t)U_{\al}x_{\pm\al_i}(-t)
|t\in\bbC,\,\,\al\in \Delta^{\rm re}_{\pm}\setminus
\{\pm\al_i\}\ran.
\]
We have the unique 
decomposition;
$U^-=U_i^-\cdot Y_{\pm\al_i}=U_{-\al_i}\cdot U^i_-.$
By using this decomposition, we get 
the canonical projection 
$\xi_i:U^-\rightarrow U_{-\al_i}$ and 
define the function $\chi_i$ on $U^-$ by 
\begin{equation}
\chi_i:=y_i^{-1}\circ\xi_i:
U^-\longrightarrow U_{-\al_i}\mapright{\sim}
\bbC,
\label{chi}
\end{equation}
and extend this to the function on 
$B^-$ by $\chi_i(u\cdot t):=\chi_i(u)$ for 
$u\in U^-$ and $t\in T$.
For a unipotent $G$-crystal $\bf(X,f_X)$, we define a function
$\vep_i:=\vep_i^X:X\rightarrow \bbC$ by 
\[
\vep_i:=\chi_i\circ{\bf f_X},
\]
and a rational function
$\gamma_i:X\rightarrow \bbC$ by 
\begin{equation}
\gamma_i:=
\al_i\circ{\rm proj}_T\circ{\bf f_X}:X\rightarrow B^-\rightarrow T\to\bbC,
\label{gammax}
\end{equation}
where ${\rm proj}_T$ is the canonical projection.\\
{\bf Remark.}
Note that the function $\vep_i$ is denoted by $\vp_i$ in
\cite{BK, N}.

Suppose that the function $\vep_i$ is not identically zero on $X$. 
We define a morphism $e_i:\bbC^\tm\tm X\rightarrow X$ by
\begin{equation}
e^c_i(x):=x_i
\left({\frac{c-1}{\vep_i(x)}}\right)(x).
\label{ei}
\end{equation}
\begin{thm}[\cite{BK},\cite{N2}]
\label{U-G}
For a unipotent $G$-crystal $({\bf X},f_X)$, 
suppose that 
the function $\vep_i$ is not identically zero
for any $i\in I$.
Then the rational functions $\gamma_i,\vep_i:X\rightarrow \bbC$ 
and 
$e_i:\bbC^\tm\tm \bf X\rightarrow X$ as above 
define a geometric 
$G$-crystal $({\bf X},\{e_i\}_{i\in I},\{\gamma_i\}_{i\in I}, 
\{\vep_i\}_{i\in I})$,
which is called the induced geometric $G$-crystals by 
unipotent $G$-crystal $({\bf X},f_X)$.
\end{thm}

\begin{pro}[\cite{BK},\cite{N2}]
For unipotent $G$-crystals $({ X},f_X)$
and $({ Y},f_Y)$, set 
the product $({Z},f_Z):=({X},f_X)
\tm({Y},f_Y)$, where 
$Z=X\tm Y$. Let $(Z,\{e^Z_i\}_{i\in I},\{\gamma^Z_i\}_{i\in I}, 
\{\vep^Z_i\}_{i\in I})$ 
be the induced geometric 
$G$-crystal from $({\bf Z},f_Z)$.
Then we obtain:
\begin{enumerate}
\item For each $i\in I$, $(x,y)\in Z$, 
\begin{equation}
\gamma^Z_i(x,y)=\gamma^X_i(x)\gamma^Y_i(y),\qq
\vep^Z_i(x,y)=
\vep^X_i(x)+\frac{\vep^Y_i(y)}{\gamma^X_i(x)}.
\label{gamma-zxy}
\end{equation}
\item
For any $i\in I$, the action 
$e^Z_i:\bbC^\tm\tm Z\rightarrow Z$ is given by: \\
$(e^Z_i)^c(x,y)=((e^X_i)^{c_1}(x),(e^Y_i)^{c_2}(y))$, where 
\begin{equation}
c_1=
\frac{c\gamma^X_i(x)\vep^X_i(x)+\vep^Y_i(y)}
{\gamma^X_i(x)\vep^X_i(x)+\vep^Y_i(y)},\,\,
c_2=
\frac{c(\gamma^X_i(x)\vep^X_i(x)+\vep^Y_i(y))}
{c\gamma^X_i(x)\vep^X_i(x)+\vep^Y_i(y)}
\label{c1c2}
\end{equation}
\end{enumerate}
\end{pro}
Here note that $c_1c_2=c$. 
The formula $c_1$ and $c_2$ in \cite{BK} 
seem to be different from ours.

\renewcommand{\thesection}{\arabic{section}}
\section{Geometric crystals on Flag variety and Schubert variety}
\label{sch}
\setcounter{equation}{0}
\renewcommand{\theequation}{\thesection.\arabic{equation}}

Let $X:=G/B$ be the flag variety, which has the cell decomposition 
$X=\sqcup_{w\in W} X_w$. Each cell $X_w$ is called a Schubert cell
associated with a Weyl group element $w\in W$. Its closure $\ovl X_w$ in $X$ is
called a Schubert variety which satisfies the closure relation
$\ovl X_w=\sqcup_{y\leq w}X_y$.
As we have seen in \cite{N},  we can associate 
geometric crystal structure with the Schubert cell
(resp. variety) $X_w$ (resp. $\ovl X_w$).

The geometric crystal on 
$X_w$ is realized in 
$B^-$ as follows:

\nd
Let $\io:={i_1}\cd {i_k}$ be one of the reduced 
expressions of $w\in W$.
Suppose that an element $w\in W$ satisfies that 
$I=I(w):=\{i_1,\cd,i_k\}$.
Define 
\[
 B^-_{\io}:=\{Y_\io(c_1,\cd,c_k):=
Y_{i_1}(c_{1})
\cd Y_{i_k}(c_{k})
\in B^-
|c_i\in \bbC^\tm\}.
\]
where $Y_i(c)=y_i(\frac{1}{c})\al^\vee_i(c)$.
The Schubert cell $X_w$ (resp. The Schubert variety $\ovl X_w$)
and $B^-_\io$ are 
birationally equivalent and 
they are isomorphic 
as induced geometric crystals.

Indeed, we describe the explicit feature of geometric 
crystal structure on $B^-_\io$:
\begin{eqnarray}
&& \gamma_i\left(Y_\io(c_1,\cd,c_k)\right)
=\al_i(\al_{i_1}^\vee(c_1)\cd \al_{i_k}^\vee(c_k))
=c_1^{{\bf a}_{i_1,i}}\cd c_k^{{\bf a}_{i_k,i}},\\
&&\vep_i(Y_{\io}(c_1,\cd,c_k))=
\sum_{1\leq j\leq k,i_j=i}
\frac{1}{c_1^{{\bf a}_{i_1,i}}\cd c_{j-1}^{{\bf a}_{i_{j-1},i}}c_j},
\label{sch-ep}\\
&&
e_i^c(Y_\io(c_1,\cd,c_k))
=:Y_\io({\mathcal C}_1,\cd,{\mathcal C}_k),\nn
\end{eqnarray}
where
\begin{equation}
{\mathcal C}_j:=
c_j\cdot \frac{\displaystyle \sum_{1\leq m\leq j,i_m=i}
 \frac{c}
{c_1^{{\bf a}_{i_1,i}}\cd c_{m-1}^{{\bf a}_{i_{m-1},i}}c_m}
+\sum_{j< m\leq k,i_m=i} \frac{1}
{c_1^{{\bf a}_{i_1,i}}\cd c_{m-1}^{{\bf a}_{i_{m-1},i}}c_m}}
{\displaystyle\sum_{1\leq m<j,i_m=i} 
 \frac{c}
{c_1^{{\bf a}_{i_1,i}}\cd c_{m-1}^{{\bf a}_{i_{m-1},i}}c_m}+
\mathop\sum_{j\leq m\leq k,i_m=i}  \frac{1}
{c_1^{{\bf a}_{i_1,i}}\cd c_{m-1}^{{\bf a}_{i_{m-1},i}}c_m}}.
\label{eici}
\end{equation}
In the case $\ge$ is semi-simple, 
we know that the flag variety $X=G/B$ coincides with 
the Schubert variety $\ovl X_{w_0}$
for the longest element $w_0$ in the Weyl group.
Thus, we have
\begin{cor}
For a semi-simple $\ge$, we have the geometric crystal
structure on the flag variety $X:=G/B$.
\end{cor}

\renewcommand{\thesection}{\arabic{section}}
\section{Geometric Crystals on $U^-$}\label{u-c}
\setcounter{equation}{0}
\renewcommand{\theequation}{\thesection.\arabic{equation}}

In this section, we associate a geometric/unipotent crystal structure
with unipotent subgroup $U^-$ of semi-simple algebraic group $G$.
In particular, for $G=SL_{n+1}(\bbC)$ we describe it explicitly.
The contents of this section is almost same as in \cite{N2}. But 
we shall see the whole setting again since we modified some definitions
and made typographical errors in the proofs of certain statements.
\subsection{$U$-variety structure on $U^-$}
In this subsection, suppose that 
$G$ is a Kac-Moody group as in Sect.2.
As mentioned in Sect.\ref{2}, Borel subgroup $B^-$ has a 
$U$-variety structure. By the similar manner, we define
$U$-variety structure on $U^-$.
As in \ref{uni-cry}, 
the multiplication map $m$ in $G$ induces an open embedding;
$m:U^-\times B\hookrightarrow G,$
then this is a birational isomorphism. 
Let us denote the inverse birational isomorphism by $h$;
\[
 h:G\longrightarrow U^-\times B.
\]
Then we define the rational maps 
$\pi^{--}:G\rightarrow U^-$ and 
$\pi^+:G\rightarrow B$ by 
$\pi^{--}:={\rm proj}_{U^-}\circ h$ 
and $\pi^+:={\rm proj}_B\circ h$.
Now we define the rational $U$-action $\al_{U^-}$ on $U^-$ by 
\[
 \al_{U^-}:=\pi^{--}\circ m:U\times U^-\longrightarrow U^-,
\]
Then we obtain 
\begin{lem}
\label{XU}
A pair ${\bf U}^-=(U^-,\al_{U^-})$ is a 
$U$-variety on a unipotent radical  $U^-\subset B^-$.
\end{lem}
\subsection{Bilinear form}
\label{bilin-form}

In this subsection, following \cite[9.4]{Kac}
we introduce the invariant 
bilinear form of finite dimensional 
modules. 
What we have introduced in \cite{N2} is subtly inexact.
So, let us reimburse it here.

Let $\lm\in P_+$ be a dominant integral weight and 
$L(\lm)$ be the associated irreducible highest weight
$\ge$-module with the fixed highest weight vector $u_\lm$. 
For $v\in L(\lm)$, define its expectation value $E(v)$ by 
\[
 v=E(v)u_{\lm}+\text{lower weight vectors}.
\]
Let $U(\ge)$ be the universal enveloping algebra of $\ge$ and 
$\what\omega:U(\ge)\to U(\ge)$ an anti-involution 
of $U(\ge)$ defined by 
$\what\omega(e_i)=f_i$, $\what\omega(f_i)=e_i$ and 
$\what\omega(h)=h$ for $i\in I$, $h\in \tt$.
Note that it is extended to an anti-involution of the group $G$
such that $\what\omega(x_i(c))=y_i(c)$, 
$\what\omega(y_i(c))=x_i(c)$ and $\what\omega(t)=t$
for $i\in I$, $c\in \bbC$ and $t\in T$. Now, we define a symmetric 
bilinear form $\lan\q,\q\ran$ on $L(\lm)$ by 
\[
 \lan u,v\ran=E(\what\omega(a)a'u_\lm),
\]
where $a,a'$ are elements in $U(\ge)$ such that 
$u=au_\lm$ and $v=a'u_\lm$.
This bilinear form satisfies
\begin{equation}
\lan gu,v\ran=\lan u,\what\omega(g)v\ran,
\label{whatomega}
\end{equation}
where $g$ is an element in $U(\ge)$ or $G$.
\subsection{Unipotent/Geometric crystal structure on $U^-$}
\label{expli-UGC}
In order to define a unipotent crystal structure on $U^-$,
let us construct a $U$-morphism 
$\cF:U^-\rightarrow B^-$.

The multiplication map $m$ in $G$ induces an open embedding;
$m:U^-\times T\times U\hookrightarrow G,$
which is a birational isomorphism. 
Thus, by the similar way as above, we obtain 
the rational map
$\pi^{0}:G\rightarrow T$.
Here note that we have 
\begin{equation}
\pi^-(x)=\pi^{--}(x)\pi^0(x)\q(x\in G).
\label{--0}
\end{equation}
Now, we give a sufficient condition for existence of 
$U$-morphism $\cF$.
\begin{lem}[\cite{N2}]
\label{suf}
Let ${\mathcal T}:U^-\rightarrow T$ be a rational map satisfying:
\begin{equation}
{\mathcal T}(\pi^{--}(xu))=\pi^0(xu){\mathcal T}(u),
\q\text{for $x\in U$ and $u\in U^-$. }
\label{T}
\end{equation}
Defining a morphism $\cF:U^-\rightarrow B^-$ by
\begin{equation}
\begin{array}{cccc}
\cF:&U^-&\longrightarrow &B^-\\
&u&\mapsto& u{\mathcal T}(u),
\end{array}
\end{equation}
then $\cF$ is a $U$-morphism $U^-\rightarrow B^-$.
\end{lem}
{\sl Proof.}
We may show 
\begin{equation}
\cF(\al_{U^-}(x,u))=\al_{B^-}(x,\cF(u)),
\q\text{for $x\in U$ and $u\in U^-$}.
\label{FU}
\end{equation}
As for the left-hand side of (\ref{FU}), we have
\[
\cF(\al_{U^-}(x,u))= \pi^{--}(xu){\mathcal T}(\pi^{--}(xu))
=\pi^{--}(xu)\pi^0(xu){\mathcal T}(u),
\]
where the last equality is due to (\ref{T}).
On the other hand, the right-hand side
 of (\ref{FU}) is written by:
\[
 \al_{B^-}(x,\cF(u))=\pi^-(xu{\mathcal T}(u))
=\pi^{--}(xu{\mathcal T}(u))\pi^0(xu{\mathcal T}(u))
=\pi^{--}(xu)\pi^0(xu){\mathcal T}(u)
\]
where the second equality is due to (\ref{--0}) and the 
third equality is obtained by the fact that 
${\mathcal T}(u)\in T\subset B$.
Now we get (\ref{FU}).\qed

Let us verify that there exists such $U$-morphism $\cF$
or rational map ${\mathcal T}$ for semi-simple cases.
Suppose that $G$ (resp. $\ge$)is semi-simple in the rest of this section.

Let $\Lm_i$ ($i=1,\cd,n$)
be a fundamental weight and $L(\Lm_i)$
be a corresponding irreducible highest weight 
$\ge$-module, where $\ge$ is a complex semi-simple
Lie algebra associated with $G$.
Let $v_\lm$ be a lowest weight vector in $L(\lm)$ such that 
$\lan v_\lm,v_\lm\ran=1$. 
Now, let us define 
a rational function $F^{(n)}_i:U^-\rightarrow \bbC$ $(i\in I)$ by 
\begin{equation}
F^{(n)}_i(u)=\lan u\cdot u_{\Lm_i},v_{\Lm_i}\ran\q
(u\in U^-).
\label{fi}
\end{equation}
We define a rational map ${\mathcal T}:U^-\rightarrow T$ by
\begin{equation}
{\mathcal T}(u):=\prod_{i\in I}\al_i^\vee(F^{(n)}_i(u)^{-1}).
\label{TT}
\end{equation}
and define a morphism $\cF:U^-\rightarrow B^-$ by
\begin{equation}
\cF(u):=u\cdot \prod_{i\in I}\al_i^\vee(F^{(n)}_i(u)^{-1}).
\label{FF}
\end{equation}
\begin{lem}
\label{cF}
The morphism $\cF:U^-\rightarrow B^-$  is a $U$-morphism.
\end{lem}
We have mentioned this statement in \cite{N2}. Nevertheless, 
since we modified the definition of the bilinear form and 
there are several typographical errors in the proof, 
we shall give a proof of this lemma again here. 

{\sl Proof of Lemma \ref{cF}.}
Let us verify that $\mathcal T$ satisfies (\ref{T}).
For $x\in U$ and $u\in U^-$ such that 
$xu\in Im(U^-\times T\times U\hookrightarrow G)$,
let $u^-\in U^-$, $u^0\in T$ and $u^+\in U$  
be the unique elements satisfying $u^-u^0u^+=xu$, 
{\it i.e.,} 
$\pi^{--}(xu)=u^-$, $\pi^0(xu)=u^0$ and $\pi(xu)=u^+$.
By (\ref{whatomega}) and the fact that $g\cdot v_{\Lm_i}=v_{\Lm_i}$ for 
any $g\in U^-$, we have 
\begin{equation}
\lan xu\cdot u_{\Lm_i}, v_{\Lm_i}\ran
=\lan u\cdot u_{\Lm_i}, \what\omega(x)\cdot v_{\Lm_i}\ran
=\lan u\cdot u_{\Lm_i}, v_{\Lm_i}\ran.
\label{xu*}
\end{equation}
On the other hand, since $g\cdot u_{\Lm_i}=u_{\Lm_i}$ for $g\in U$,
we have 
\begin{equation}
\begin{array}{l}
 \lan xu\cdot u_{\Lm_i}, v_{\Lm_i}\ran
=\lan \pi^{--}(xu)\pi^0(xu)\pi(xu)\cdot u_{\Lm_i}, v_{\Lm_i}\ran
\\
\qq=\lan \pi^{--}(xu)\pi^0(xu)\cdot u_{\Lm_i}, v_{\Lm_i}\ran
=\Lm_i(\pi^0(xu))\lan \pi^{--}(xu)\cdot u_{\Lm_i}, v_{\Lm_i}\ran,
\end{array}
\label{xu0}
\end{equation}
where we regard $\Lm_i$ as an element in $X^*(T)$ such that 
$\Lm_i(\al_j^\vee(c))=c^{\del_{i,j}}$.
Hence, by (\ref{xu*}), (\ref{xu0}), we have
\begin{eqnarray*}
F^{(n)}_i(\pi^{--}(xu))&=&\lan \pi^{--}(xu)\cdot u_{\Lm_i}, v_{\Lm_i}\ran
=\Lm_i(\pi^0(xu))^{-1}\lan xu\cdot u_{\Lm_i}, v_{\Lm_i}\ran\\
&=&\Lm_i(\pi^0(xu))^{-1}\lan u\cdot u_{\Lm_i}, v_{\Lm_i}\ran
=\Lm_i(\pi^0(xu))^{-1}F^{(n)}_i(u).
\end{eqnarray*}
By the formula
\[
 \prod_i\al_i^\vee(\Lm_i(t))=t, \q(t\in T),
\]
and the definitions of $\mathcal T$ and $\cF$, 
we obtained (\ref{T}).\qed

\begin{cor}
Suppose that $G$ $($resp. $\ge)$ is semi-simple. 
Then $(U^-,F)$ is a unipotent crystal.
\end{cor}
As we have seen in \ref{u->g}, we can associate 
geometric crystal structure with the unipotent subgroup $U^-$
since it has a unipotent crystal structure.

Let us denote the function $\chi_i:U^-\to\bbC$ in (\ref{chi}) by 
$\vep_i:U^-\to\bbC$ here.
It is trivial that the function 
$\vep_i:U^-\rightarrow \bbC$ is not identically zero.
Thus, defining the morphisms 
$e_i:\bbC^\times\times U^-\rightarrow U^-$ and 
$\gamma_i:U^-\rightarrow \bbC$ by
\begin{equation}
e_i(c,u)=e_i^c(u):=x_i(\frac{c-1}{\vep_i(u)})(u),\qq
\gamma_i(u):=\al_i({\mathcal T}(u)),\qq
\text{($u\in U^-$ and $c\in \bbC^\times$)},
\end{equation}
It follows from Theorem \ref{U-G}:
\begin{thm}
If $G$ is semi-simple, then 
$\chi_{U^-}:=(U^-,\{\gamma_i\}_{i\in I},\{\vep_i\}_{i\in I},
\{e_i\}_{i\in I})$ is a geometric 
crystal.
\end{thm}

\subsection{Explicit Form of the Geometric Crystal Structure of $U^-$}
\label{subsec-ex-uc}

Let $\io_0=i_1i_2\cd i_N$ be a reduced longest word of 
a semi-simple Lie algebra $\ge$ and set 
\[
 U_{\io_0}^-:=\{y_{\io_0}(a)=y_{i_1}(a_1)\cd
y_{i_k}(a_N)|a_1,\cd, a_N\in \bbC\},
\]
which is birationally equivalent to $U^-$.
Thus, using this we describe an
explicit form of the geometric crystal structure of $U^-$:
For $\io_0$ and $i\in I$, define $\{j_1,j_2,\cd,j_l\}
:=\{j|1\leq j\leq k,\,\,m_j=i\}$, where 
$1\leq j_1<\cd <j_l\leq N$ and set 
\[
{L^{(i)}_m(a;c)}:=\frac{c(a_{j_1}+\cd +a_{j_m})
+a_{j_{m+1}}+\cd +a_{j_l}}
{a_{j_1}+\cd +a_{j_l}}\,\,\q(1\leq m\leq l,\,\,c\in\bbC).
\]
Then, we have
\begin{eqnarray}
&&\qq
\vep_i(y_{\io_0}(a))=\sum_{i_j=i}a_j,\qq
\gamma_i(y_{\io_0}(a))
=\al_i(\prod_j\al_j^\vee(F^{(n)}_j(y_{\io_0}(a)))^{-1}),\nn \\
&& \qq e^c_i(y_{\io_0}(a))=
x_i(\frac{c-1}{\vep_i(y_{\io_0}(a))})(y_{\io_0}(a))
=y_{i_1}({a'_1})\cd y_{i_N}
({a'_N}),\label{ex-uc}\\
&&\text{where\q}
{a'_{j_m}}=\frac{a_{j_m}}
{L^{(i)}_{m-1}(a;c)L^{(i)}_m(a;c)}\q
(i_{j_m}=i),
\q{a'_{i_p}}=\frac{a_{p}}
{{L^{(i)}_{m-1}(a;c)}^{{\bf a}_{i,i_p}}}\q
(j_{m-1}<p< j_m).\nn
\end{eqnarray}
Note that $L_0^{(i)}(a;c)=1$.

\renewcommand{\thesection}{\arabic{section}}
\section{Fundamental Representations}
\setcounter{equation}{0}
\renewcommand{\theequation}{\thesection.\arabic{equation}}

In order to get the explicit form of the function $F^{(n)}_i$
in the next section, 
we shall see some technical lemmas in this section.

\subsection{Type $C_n$}

Let $I:=\{1,2,\cd\,n\}$ be the index set of the 
simple roots of type $C_n$. The Cartan matrix 
$A=({\bf a}_{i,j})_{i,j\in I}$ of type 
$C_n$ is given by 
\[
 {\bf a}_{i,j}=\begin{cases}
2&\text{if }i=j,\\
-1&\text{if }|i-j|=1\text{ and }(i,j)\ne(n-1,n)\\
-2&\text{if }(i,j)=(n-1,n),\\
0&\text{otherwise}.
\end{cases}
\]
Here $\al_i$ ($i\ne n$) is a short root and 
$\al_n$ is the long root.
Let $\{h_i\}_{i\in I}$ be the set of simple co-roots
and $\{\Lm_i\}_{i\in I}$ be the set of fundamental 
weights satisfying $\al_j(h_i)={\bf a}_{i,j}$ 
and $\Lm_i(h_j)=\del_{i,j}$.

First, let us describe the vector representation 
$L(\Lm_1)$. Set ${\mathbf B}^{(n)}:=
\{v_j,v_{\ovl j}|j=1,2,\cd,n.\}$. The weight of $v_j$ is as follows:
\[
 {\rm wt}(v_j)=\begin{cases}
\Lm_i-\Lm_{i-1}&\text{ if }i=1,\cd,n,\\
\Lm_{i-1}-\Lm_{i}&\text{ if }i=\ovl 1,\cd,\ovl n,
\end{cases}
\]
where $\Lm_0=0$.
The actions of $e_i$ and $f_i$ 
are given by:
\begin{eqnarray}
&&f_iv_i=v_{i+1},\q f_iv_{\ovl i+1}=v_{\ovl i},\q
e_iv_{i+1}=v_i,\q e_iv_{\ovl i}=v_{\ovl{i+1}}
\q(1\leq i<n),\label{c-f1}\\
&&f_nv_n=v_{\ovl n},\qq 
e_nv_{\ovl n}=v_n,\label{c-f2}
\end{eqnarray}
and the other actions are trivial.

Let $\Lm_i^{(n)}$ be the $i$-th fundamental weight of type 
$C_n$, where we add the superscript $(n)$ to emphasize the
rank of the corresponding Lie algebra.
As is well-known that the fundamental representation 
$L(\Lm_i^{(n)})$ $(1\leq i\leq n)$
is embedded in $L(\Lm_1^{(n)})^{\ot i}$
with multiplicity free.
The explicit form of the highest(resp. lowest) weight 
vector $u_{\Lm_i^{(n)}}$ (resp. $v_{\Lm_i^{(n)}}$)
of $L(\Lm_i^{(n)}$ is realized in 
$L(\Lm_1^{(n)})^{\ot i}$ as follows:
\begin{equation}
\begin{array}{ccc}\displaystyle
u_{\Lm_i^{(n)}}&=&\displaystyle
\sum_{\sigma\in{\mathfrak S}_i}{\rm sgn}(\sigma)
v_{\sigma(1)}\ot\cd\ot v_{\sigma(i)},\\
v_{\Lm_i^{(n)}}&=&\displaystyle
\sum_{\sigma\in{\mathfrak S}_i}{\rm sgn}(\sigma)
v_{\ovl{\sigma(i)}}\ot\cd\ot v_{\ovl{\sigma(1)}}, 
\end{array}
\label{h-l}
\end{equation}
where ${\mathfrak S}_i$ is the $i$-th symmetric group.
For $x=\sum_{i_1,\cd,i_k}c_{i_1,\cd,i_k}v_{i_1}
\ot\cd\ot v_{i_k}
\in L(\Lm_1^{(n)})^{\ot k}$, $v\in L(\Lm_1^{(n)})$ and
$j\in\{1,\cd,k\}$,
let us define:
\[
 x[v;j]:=\sum_{i_1,\cd,i_k}c_{i_1,\cd,i_k}v_{i_1}\ot\cd
v_{i_{j-1}}\ot v\ot v_{i_j}\ot\cd\ot v_{i_k}\in 
L(\Lm_1^{(n)})^{\ot k+1}.
\]
Let $u'$ (resp. $v'$) be the vector in 
$L(\Lm_1^{(n+1)})^{\ot i}$ $(i<n)$ whose 
explicit form is given by replacing 
$v_j$ (resp. $v_{\ovl j}$)$\in L(\Lm_1^{(n)})$ by $v_{j+1}$ 
(resp. $v_{\ovl{j+1}}$) $\in L(\Lm_1^{(n+1)})$
($j=1,\cd,i$) in the vector
$u_{\Lm_i^{(n)}}$ (resp. $v_{\Lm_i^{(n)}}$) in (\ref{h-l}).
Then, they satisfy $e_i u'=0$ (resp. $f_i v'=0$) for $i=2,\cd,n$.
Here for the vectors $v_j,v_{\ovl j}$ ($i=1,\cd,n$) in 
${\mathbf B}^{(n)}$ and ${\mathbf B}^{(n+1)}$ we shall 
use the same notations. 
\begin{lem}
\label{c+}
We have
\begin{eqnarray}
u_{\Lm_{i+1}^{(n+1)}}&=&\sum_{j=1}^{i+1}(-1)^{j-1}u'[v_1;j],\\
v_{\Lm_{i+1}^{(n+1)}}&=&
\sum_{j=1}^{i+1}(-1)^{i-j}v'[v_{\ovl 1};j].
\end{eqnarray}
\end{lem}
{\sl Proof.}
Suppose $\sigma\in{\mathfrak S}_{i+1}$ is in the form 
\[
 \sigma=\begin{pmatrix}1&\cd&j\cd&i+1\\
\sigma(1)&\cd&1\cd&\sigma(i+1)\end{pmatrix},
\]
that is, $\sigma(j)=1$. Then we have 
\[
 \sigma=(\sigma(1),\cd,\sigma(j-1),1)
\begin{pmatrix}1&2\cd&j&j+1&\cd&i+1\\
1&\sigma(1)\cd&\sigma(j-1)&\sigma(j+1)&
\cd&\sigma(i+1)\end{pmatrix},
\]
where $(\sigma(1),\cd,\sigma(j-1),1)$ is a cycle.
Since $\sigma(k)\ne1$ for $k\ne j$, we have that 
\[
\sigma':= \begin{pmatrix}2\cd&j&j+1\cd &i+1\\
\sigma(1)\cd&\sigma(j-1)&\sigma(j+1)
\cd&\sigma(i+1)\end{pmatrix},
\]
is a permutation of $\{2,3,\cd,i+1\}$ satisfying 
$\sigma'(k)=\sigma(k-1)$ for $2\leq k\leq j$ and 
$\sigma'(k)=\sigma(k)$ for $k>j$.
Thus, we have ${\rm sgn}(\sigma)
=(-1)^{j-1}{\rm sgn}(\sigma')$.
Hence, we have
\begin{eqnarray*}
u_{\Lm_{i+1}^{(n+1)}}&=&
\sum_{j=1}^{i+1}\sum_{\sigma'\in{\mathfrak S}'}
(-1)^{j-1}{\rm sgn}(\sigma')v_{\sigma'(2)}
\ot\cd\ot v_{\sigma'(j)}
\ot v_1\ot v_{\sigma'(j+1)}\ot \cd\ot v_{\sigma'(i+1)},\\
&=&\sum_{j=1}^{i+1}(-1)^{j-1}u'[v_1;j],
\end{eqnarray*}
where ${\mathfrak S}'={\mathfrak S}_{\{2,\cd,i+1\}}$. 
The case of $v_{\Lm_{i+1}^{(n+1)}}$ is shown similarly.\qed
\subsection{Type $B_n$}\label{Bn}

Let $I:=\{1,2,\cd\,n\}$ be the index set of the 
simple roots of type $B_n$. The Cartan matrix 
$A=({\bf a}_{i,j})_{i,j\in I}$ of type 
$B_n$ is given by 
\[
 {\bf a}_{i,j}=\begin{cases}
2&\text{if }i=j,\\
-1&\text{if }|i-j|=1\text{ and }(i,j)\ne(n,n-1)\\
-2&\text{if }(i,j)=(n,n-1),\\
0&\text{otherwise}.
\end{cases}
\]
Here $\al_i$ ($i\ne n$) is a long root and 
$\al_n$ is the short root.
Let $\{h_i\}_{i\in I}$ be the set of simple co-roots
and $\{\Lm_i\}_{i\in I}$ be the set of fundamental 
weights satisfying $\al_j(h_i)={\bf a}_{i,j}$ 
and $\Lm_i(h_j)=\del_{i,j}$.

First, let us describe the vector representation $L(\Lm_1)$ for 
$B_n$.
Set ${\mathbf B}^{(n)}:=\{v_j,\,v_{\ovl j}|j=1,2,\cd,n\}\cup\{v_0\}$. 
The weight of $v_j$ is as follows:
\begin{eqnarray*}
&& {\rm wt}(v_j)=
\Lm_i-\Lm_{i-1},\q
{\rm wt}(v_{\ovl j})=\Lm_{i-1}-\Lm_{i}\q
(i=1,\cd,n-1),\\
&&{\rm wt}(v_n)=
2\Lm_n-\Lm_{n-1},\q
{\rm wt}(v_{\ovl n})=
\Lm_{n-1}-2\Lm_n,\q {\rm wt}(v_0)=0,
\end{eqnarray*}
where $\Lm_0=0$.
The actions of $e_i$ and $f_i$ 
are given by:
\begin{eqnarray}
&&f_iv_i=v_{i+1},\q f_iv_{\ovl{i+1}}=v_{\ovl i},\q
e_iv_{i+1}=v_i,\q e_iv_{\ovl i}=v_{\ovl{i+1}}
\q(1\leq i<n),\label{b-f1}\\
&&f_nv_n=v_0,\q f_nv_0=2v_{\ovl n},\q 
e_nv_0=2v_n,\q e_nv_{\ovl n}=v_0,\label{b-f2}
\end{eqnarray}
and the other actions are trivial.

The last fundamental representation $L(\Lm_n)$ 
is called the ``spin representation'' whose dimension
is $2^n$. It is realized as follows:
Set $V^{(n)}_{sp}:=\bigoplus_{\ep\in B^{(n)}_{sp}}\bbC\ep$ and 
\[
 B^{(n)}_{sp}:=\{(\ep_1,\cd,\ep_n)|\ep_i\in\{+,-\}\}.
\]
Define the explicit action of $h_i$, $e_i$ and $f_i$ on 
$V^{(n)}_{sp}$ by 
\begin{eqnarray}
h_i(\ep_1,\cd,\ep_n)&=&\begin{cases}
\frac{\ep_i\cdot1-\ep_{i+1}\cdot1}{2}
(\ep_1,\cd,\ep_n),&\text{ if }i\ne n,\\
\ep_n(\ep_1,\cd,\ep_n)&\text{ if }i=n,
\end{cases}\\
f_i(\ep_1,\cd,\ep_n)&=&\begin{cases}
(\cd,\buildrel{i}\over-,\buildrel{i+1}\over+,\cd)
&\text{ if }\ep_i=+,\,\,\ep_{i+1}=-,\,\,
i\ne n,\\
(\cd\cd\cd,\buildrel{n}\over-)
&\text{ if }\ep_n=+,\,\,i=n,\\
\qq0&\text{otherwise}
\end{cases}\\
e_i(\ep_1,\cd,\ep_n)&=&\begin{cases}
(\cd,\buildrel{i}\over+,\buildrel{i+1}\over-,\cd)
&\text{ if }\ep_i=-,\,\,\ep_{i+1}=+,\,\,
i\ne n,\\
(\cd\cd\cd,\buildrel{n}\over+)
&\text{ if }\ep_n=-,\,\,i=n,\\
\qq0&\text{otherwise.}
\end{cases}
\end{eqnarray}
Then the module $V^{(n)}_{sp}$ is isomorphic to 
$L(\Lm_n)$ as a $B_n$-module.

The following decomposition is well-known:
\[
 L(\Lm_n)\ot L(\Lm_n)\cong
L(0)\oplus L(\Lm_1)\oplus L(\Lm_2)
\oplus\cd\oplus L(\Lm_{n-1})
\oplus L(2\Lm_n).
\]
Let us describe the explicit form of the 
highest (resp. lowest)weight vector 
$u_i^{(n)}$(resp. $v_i^{(n)}$) in $L(\Lm_i^{(n)})$
($i=1,\cd,n-1$) by using the vectors in $B^{(n)}_{sp}$,
where $\Lm_i^{(n)}=\Lm_i$ and we emphasize the rank 
of the Lie algebra by adding the superscript $(n)$.
For $\ep=(\ep_1,\cd,\ep_n)\in B^{(n)}_{sp}$, we define the 
signature sg($\ep$) as follows:
set $J(\ep):=\{j_1,\cd,j_m\}\subset\{1,\cd,n\}$ such that 
$\ep_{j_k}=-$ for any $k=1,\cd,m$ and if $p\ne j_k$, 
then $\ep_p=+$.
Then we define 
\[
 {\rm sg}(\ep)=(-1)^{\sum_{k=1}^{m}(n-j_k+1)}.
\]
It is easy to see:
\begin{equation}
u_i^{(n)}:=\sum_{\ep,\ep'\text{ satisfies (H)}}
{\rm sg}(\ep)\ep\ot\ep',\qq
v_i^{(n)}:=\sum_{\ep,\ep'\text{ satisfies (L)}}
(-1)^{\frac{n(n+1)}{2}}{\rm sg}(\ep')\ep\ot\ep',
\end{equation}
where for $\ep=(\ep_1,\cd,\ep_n)$ and $\ep'=(\ep'_1,\cd,\ep'_n)$, 
\begin{enumerate}
\item[(H)] $\ep_j=\ep'_j=+$
($j\leq i$) and $\ep_j\cdot\ep'_j=-$ ($j>i$).
\item[(L)]
$\ep_j=\ep'_j=-$
($j\leq i$) and $\ep_j\cdot\ep'_j=-$ ($j>i$).
\end{enumerate}
For $v=\sum_{\ep,\ep'\in B^{(n)}_{sp}}c_{\ep,\ep'}\ep\ot\ep' \in
V_{sp}^{(n)}\ot V_{sp}^{(n)}$, 
let us denote 
\[
 (\pm,v):=\sum_{\ep,\ep'\in B^{(n)}_{sp}}c_{\ep,\ep'}(\pm,\ep)
\ot(\pm,\ep')\in V_{sp}^{(n+1)}\ot V_{sp}^{(n+1)},
\]
where $(\pm,\ep)\in B^{(n+1)}_{sp}$.

\nd
{\sl Remark.} Note that for $B^{(n)}_{sp}$ and 
$i=2,\cd,n+1$, we have $e_i(\pm,v)=(\pm,e_{i-1}v)$
and $f_i(\pm,v)=(\pm,f_{i-1}v)$.
\begin{lem}
\label{b+-}
We have
\begin{equation}
u_{i+1}^{(n+1)}=(+,u_i^{(n)}),\qq
v_{i+1}^{(n+1)}=(-,v_i^{(n)})
\q(i=1,\cd,n).
\end{equation}
\end{lem}
{\sl Proof.} 
For $\ep=(\ep_1,\cd,\ep_n)\in B^{(n)}_{sp}$, 
let $J(\ep)=\{j_1,\cd,j_m\}$ be the same set as above.
Then $J(+,\ep)=\{j_1+1,j_2+1,\cd,j_m+1\}$.
Thus, we have $u_{i+1}^{(n+1)}=(+,u_i^{(n)})$ since 
\[
 {\rm sg}(+,\ep)=(-1)^{\sum_{k=1}^{m}((n+1)-(j_k+1)+1)}
=(-1)^{\sum_{k=1}^{m}(n-j_k+1)}={\rm sg}(\ep).
\]
Next, let us show the case $v_j^{(n)}$. 
For $\ep=(\ep_1,\cd,\ep_n)\in B^{(n)}_{sp}$, 
let $J(\ep)=\{j_1,\cd,j_m\}$. 
Then we get 
$J(-,\ep)=\{1,j_1+1,j_2+1,\cd,j_m+1\}$.
Thus, we obtain
\[
 {\rm sg}(-,\ep)=(-1)^{((n+1)-1+1)
+\sum_{k=1}^{m}((n+1)-(j_k+1)+1)}
=(-1)^{n+1}\cdot(-1)^{\sum_{k=1}^{m}n-j_k+1}
=(-1)^{n+1}{\rm sg}(\ep),
\]
which implies $(-,v_i^{(n)})=v_{i+1}^{(n+1)}$.\qed

\subsection{Type $D_n$}\label{Dn}

Let $I:=\{1,2,\cd\,n\}$ be the index set of the 
simple roots of type $D_n$. The Cartan matrix 
$A=({\bf a}_{i,j})_{i,j\in I}$ of type 
$D_n$ is as follows:
\[
 {\bf a}_{i,j}=\begin{cases}
2&\text{if }i=j,\\
-1&\text{if }|i-j|=1
\text{ and }(i,j)\ne(n,n-1),\,(n-1,n), 
\text{ or }(i,j)=(n-2,n),\,\,(n,n-2)\\
0&\text{otherwise}.
\end{cases}
\]
Let $\{h_i\}_{i\in I}$ be the set of simple co-roots
and $\{\Lm_i\}_{i\in I}$ be the set of fundamental 
weights satisfying $\al_j(h_i)={\bf a}_{i,j}$ 
and $\Lm_i(h_j)=\del_{i,j}$.

First, let us describe the vector representation $L(\Lm_1)$ for 
$D_n$.
Set ${\mathbf B}^{(n)}:=\{v_j, v_{\ovl j}|j=1,2,\cd,n\}$. 
The weight of $v_j$ is as follows:
\begin{eqnarray*}
&& {\rm wt}(v_j)=
\Lm_i-\Lm_{i-1},\q
{\rm wt}(v_{\ovl j})=\Lm_{i-1}-\Lm_{i}\q
(i=1,\cd,n-1),\\
&&{\rm wt}(v_n)=
\Lm_{n-1}+\Lm_n-\Lm_{n-2},\q
{\rm wt}(v_{\ovl n})=
\Lm_{n-2}-\Lm_{n-1}+\Lm_n,
\end{eqnarray*}
where $\Lm_0=0$.
The actions of $e_i$ and $f_i$ 
are given by:
\begin{eqnarray}
&&f_iv_i=v_{i+1},\q f_iv_{\ovl{i+1}}=v_{\ovl i},\q
e_iv_{i+1}=v_i,\q e_iv_{\ovl i}=v_{\ovl{i+1}}
\q(1\leq i<n),\\
&&f_nv_n=v_{\ovl{n-1}},\q f_{n-1}v_{\ovl n}=v_{\ovl{n-1}},\q 
e_{n-1}v_{\ovl {n-1}}=v_{\ovl n},\q e_nv_{\ovl{n-1}}=v_n,
\end{eqnarray}
and the other actions are trivial.

The last two fundamental representations $L(\Lm_{n-1})$ and $L(\Lm_n)$ 
are also called the ``spin representations'' whose dimension
are $2^{n-1}$. They are realized as follows:
Set $V^{(+,n)}_{sp}$(resp. $V^{(-,n)}_{sp}$)
$:=\bigoplus_{\ep\in B^{(+,n)}_{sp}(\text{resp. }B^{(-,n)}_{sp})}\bbC\ep$ and 
\[
 B^{(+,n)}_{sp}(\text{resp. }B^{(-,n)}_{sp})
:=\{(\ep_1,\cd,\ep_n)|\ep_i\in\{+,-\},\ep_1\cd\ep_n=+(\text{resp. }-)\}.
\]
Define the explicit action of $h_i$, $e_i$ and $f_i$ on 
$V^{(\pm,n)}_{sp}$ by 
\begin{eqnarray}
h_i(\ep_1,\cd,\ep_n)&=&\begin{cases}
\frac{\ep_i\cdot1-\ep_{i+1}\cdot1}{2}
(\ep_1,\cd,\ep_n),&\text{ if }i\ne n,\\
\frac{\ep_{n-1}\cdot1+\ep_n\cdot1}{2}(\ep_1,\cd,\ep_n)&\text{ if }i=n,
\end{cases}\\
f_i(\ep_1,\cd,\ep_n)&=&\begin{cases}
(\cd,\buildrel{i}\over-,\buildrel{i+1}\over+,\cd)
&\text{ if }\ep_i=+,\,\,\ep_{i+1}=-,\,\,
i\ne n,\\
(\cd\cd\cd,\buildrel{n-1}\over-,\buildrel{n}\over-)
&\text{ if }\ep_{n-1}=+,\ep_n=+,\,\,i=n,\\
\qq0&\text{otherwise}
\end{cases}\\
e_i(\ep_1,\cd,\ep_n)&=&\begin{cases}
(\cd,\buildrel{i}\over+,\buildrel{i+1}\over-,\cd)
&\text{ if }\ep_i=-,\,\,\ep_{i+1}=+,\,\,
i\ne n,\\
(\cd\cd\cd,\buildrel{n-1}\over+,\buildrel{n}\over+)
&\text{ if }\ep_{n-1}=-,\ep_n=-,\,\,i=n,\\
\qq0&\text{otherwise.}
\end{cases}
\end{eqnarray}
Then the module $V^{(+,n)}_{sp}$ (resp. $V^{(-,n)}_{sp}$) is isomorphic to 
$L(\Lm_n)$ (resp. $L(\Lm_{n-1})$) as a $D_n$-module.

The following decomposition is well-known:
\begin{eqnarray}
&&V^{(+,n)}_{sp}\ot V^{(+,n)}_{sp}\cong
\begin{cases}L(0)\oplus L(\Lm_2)
\oplus\cd\oplus L(\Lm_{2m-2})\oplus L(2\lm_{2m})
&\text{ if } n=2m,\\
L(\Lm_1)\oplus L(\Lm_3)
\oplus\cd\oplus L(\Lm_{2m-1})\oplus L(2\lm_{2m+1})
&\text{ if } n=2m+1,
\end{cases}
\\
&&
V^{(+,n)}_{sp}\ot V^{(-,n)}_{sp}\cong
\begin{cases}L(\Lm_1)\oplus L(\Lm_3)
\oplus\cd\oplus L(\Lm_{2m-3})\oplus L(2\lm_{2m-1})
&\text{ if } n=2m,\\
L(0)\oplus L(\Lm_2)
\oplus\cd\oplus L(\Lm_{2m-2})\oplus L(2\lm_{2m})
&\text{ if } n=2m+1,
\end{cases}
\end{eqnarray}
Thus, we know that each fundamental representation 
$L(\Lm_i)$ ($i=1,\cd,n-2$) is embedded in 
$V^{(+,n)}_{sp}\ot V^{(\pm,n)}_{sp}$ with multiplicity free.
Now, let us describe the unique (up to constant) highest 
(resp. lowest) weight vector $u_i^{(n)}$ (resp. $v_i^{(n)}$)
in $L(\Lm_i)$.
It is trivial that 
\begin{eqnarray*}
&& u^{(n)}_n=(+,\cd,+,+),\q
 v^{(n)}_n=\begin{cases}(-,\cd,-,-)&\text{ if }n
\text{ is even.}\\
(-,\cd,-,+)&\text{ if }n
\text{ is odd.}
\end{cases}
\\
&&u^{(n)}_{n-1}=(+,\cd,+,-),\q
 v^{(n)}_{n-1}=\begin{cases}(-,\cd,-,+)&\text{ if }n
\text{ is even.}\\
(-,\cd,-,-)&\text{ if }n
\text{ is odd.}
\end{cases}
\end{eqnarray*}

For $\ep\in V^{(\pm,n)}_{sp}$, let ${\rm sg}(\ep)$ be the same 
as above. Then it is immediate that 
\begin{equation}
u_i^{(n)}:=\sum_{\ep,\ep'\text{ satisfies (H)}}
{\rm sg}(\ep)\ep\ot\ep',\qq
v_i^{(n)}:=\sum_{\ep,\ep'\text{ satisfies (L)}}
(-1)^{\frac{n(n+1)}{2}}{\rm sg}(\ep')\ep\ot\ep',
\end{equation}
where
if $n-i$ is even (resp. odd), then $u_i^{(n)}$ and  $v_i^{(n)}$
are in $V^{(+,n)}_{sp}\ot V^{(+,n)}_{sp}$
(resp. $V^{(+,n)}_{sp}\ot V^{(-,n)}_{sp}$), and 
\begin{enumerate}
\item[(H)] $\ep_j=\ep'_j=+$
($j\leq i$) and $\ep_j\cdot\ep'_j=-$ ($j>i$).
\item[(L)]
$\ep_j=\ep'_j=-$
($j\leq i$) and $\ep_j\cdot\ep'_j=-$ ($j>i$).
\end{enumerate}
By arguing similarly to the $B_n$-case, we have
\begin{lem}
\begin{eqnarray}
&&u_{i+1}^{(n+1)}=(+,u_i^{(n)}),\qq
v_{i+1}^{(n+1)}=(-,v_i^{(n)})
\q(i=1,\cd,n-2),\\
&&u^{(n+1)}_{n}=(+,u^{(n)}_{n-1}),\q 
v^{(n+1)}_{n}=(-,v^{(n)}_{n-1}),\\
&&u^{(n+1)}_{n+1}=(+,u^{(n)}_{n}),\q 
v^{(n+1)}_{n+1}=(-,v^{(n)}_{n}),
\end{eqnarray}
where the notation $(\pm,u)$ is the same one as 
in the previous subsection.
\end{lem}
{\sl Remark.} Similar to \ref{Bn}, 
for $v\in B^{(\pm,n)}_{sp}$ and 
$i=2,\cd,n+1$, we have $e_i(\pm,v)=(\pm,e_{i-1}v)$
and $f_i(\pm,v)=(\pm,f_{i-1}v)$.

\renewcommand{\thesection}{\arabic{section}}
\section{Function $F^{(n)}_i$}
\setcounter{equation}{0}
\renewcommand{\theequation}{\thesection.\arabic{equation}}

Fix the following reduced longest word:
\[
 \io_0=\begin{cases}
(n\,\,n-1\cd 21)(n\,\,n-1\cd 32)\cd(n\,n-1)(n)&A_n,\\
(12\cd n-1\,n\,n-1\cd 21)(2\cd n-1\,n\,n-1\cd 2)\cd(n-1\,n\,n-1)(n)&
B_n,\,\,C_n\\
(12\cd n-1\,n\,n-2\cd 21)(2\cd n-1\,\,n\,\,n-2\cd 2)\cd &D_n \\
\cd(n-2\,\,n-1\,\,n\,\,n-2)(n-1\,\,n).
\end{cases}
\]
For these words and $y\in U^-_{\io_0}$, 
 we shall obtain the explicit forms of $F^{(n)}_i(y)$.
\begin{pro}
\label{Fi}
We have: \hfill\break
\vspace{-4mm}
\begin{enumerate}
\item
$A_n$-case:
For $y=(y_n(a_{1,n})\cd y_1(a_{1,1}))\cd 
(y_{n}(a_{n-1,n})y_{n-1}(a_{n-1,n}))y_n(a_{n,n}))$
\[
F^{(n)}_i(y)=\prod_{k=1}^{i}\prod_{j=k}^{n-i+k}
a_{k,j}
\]
\item
$B_n$-case:
For $y=(y_1(a_{1,1})\cd y_n(a_{1,n})
y_{n-1}(\ovl a_{1,n-1})\cd y_1(\ovl a_{1,1}))$\\
$\times 
(y_2(a_{2,2})\cd y_2(\ovl a_{2,2}))\times\cd\times$
$(y_{n-1}(a_{n-1,n-1})y_n(a_{n-1,n})
y_{n-1}(\ovl a_{n-1,n-1}))$
$y_n(a_{n,n})$
\[
F^{(n)}_i(y)=
\begin{cases}\displaystyle
\prod_{1\leq k\leq i\leq j\leq n}
a_{k,j}\ovl a_{k,j}\prod_{1\leq k\leq j<i}
\ovl a_{k,j}^2\,
&\text{if }i<n,\\
\displaystyle
\prod_{1\leq k\leq j\leq n}
\ovl a_{k,j}&\text{if }i=n,
\end{cases}
\]
where we understand $\ovl a_{k,n}=a_{k,n}$.
\item
$C_n$-case: For the same $y$ in (ii),
\[
F^{(n)}_i(y)=
\prod_{1\leq k\leq i\leq j\leq n}
a_{k,j}
\prod_{1\leq k\leq i\leq j< n}\ovl a_{k,j}
\prod_{1\leq k\leq j< i}
\ovl a_{k,j}^2,
\]
where if $i=n$, we understand the second factor is equal to $1$.
\item 
$D_n$-case:
For $y= (y_1(a_{1,1})\cd y_{n-1}(a_{1,n-1})
y_{n}(\ovl a_{1,n-1})y_{n-2}(\ovl a_{1,n-2})\cd y_1(\ovl a_{1,1}))$\\
$\qq\qq\qq\cd\cd\cd$\\
$\times 
y_{n-2}(a_{n-2,n-2})y_{n-1}(a_{n-2,n-1})
y_n(\ovl a_{n-2,n-1})y_{n-2}(\ovl a_{n-2,n-2})$\\
$\times y_{n-1}(a_{n-1,n-1})y_n(\ovl a_{n-1,n-1})$,
\[
F^{(n)}_i(y)=
\begin{cases}\displaystyle
\prod_{1\leq k\leq i\leq j<n}
a_{k,j}\ovl a_{k,j}\prod_{1\leq k\leq j< i}
\ovl a_{k,j}^2\,
&\text{if }1\leq i\leq n-2,\\
\displaystyle
(\prod_{1\leq k\leq j<n-1}
\ovl a_{k,j})\cdot(a_{1,n-1}\ovl a_{2,n-1}\cd a_{n-1,n-1})
&\text{if }i=n-1,\q n:\text{ even},\\
\displaystyle
(\prod_{1\leq k\leq j<n-1}
\ovl a_{k,j}\cdot)(\ovl a_{1,n-1}a_{2,n-1}\cd a_{n-1,n-1})
&\text{if }i=n-1,\q n:\text{ odd},\\
\displaystyle
(\prod_{1\leq k\leq j<n-1}
\ovl a_{k,j})\cdot(a_{1,n-1}\ovl a_{2,n-1}\cd\ovl a_{n-1,n-1})
&\text{if }i=n,\q n:\text{ odd},\\
\displaystyle
(\prod_{1\leq k\leq j<n-1}
\ovl a_{k,j})\cdot(\ovl a_{1,n-1} a_{2,n-1}\cd\ovl a_{n-1,n-1})
&\text{if }i=n,\q n:\text{ even},\\
\end{cases}
\]
\end{enumerate}
\end{pro}
Note that the notation $\ovl a_{i,j}$ does NOT mean the 
complex conjugate of $a_{i,j}$ but means a variable without 
any relation to $a_{i,j}$.

\vskip3mm
{\sl Proof.}
Indeed, the $A_n$-case has already been given in \cite{N2}. 
Let us see the case $i=1$ for other three types.

Let $v_1$ be the highest weight vector in $L(\Lm_1)$.
Since $f_iv_1=0$ for $i\ne 1$, 
we have 
\[
yv_1=\begin{cases}y_1(a_{1,1})\cd y_n(a_{1,n})
y_{n-1}(\ovl a_{1,n-1})\cd y_1(\ovl a_{1,1})v_1,
&B_n,\,\,C_n,\\
y_1(a_{1,1})\cd y_{n-1}(a_{1,n-1})
y_{n}(\ovl a_{1,n-1})y_{n-2}(\ovl a_{1,n-2})\cd y_1(\ovl a_{1,1})v_1,
&D_n
\end{cases}
\]
where $y\in U^-_{\io_0}$ is as in Proposition \ref{Fi}.
Since  $f_i^2=0$ on $L(\Lm_1)$ for types $C_n$ and $D_n$ and 
$f_i^2=0$ ($i\ne n$) and $f_n^3=0$ on $L(\Lm_1)$ for types $B_n$,
we obtain 
\begin{equation}
yv_1=\begin{cases}
a_{1,1}\cd a_{1,n}^2\ovl a_{1,n-1}\cd \ovl a_{1,1}v_{\ovl 1}
+w&B_n,\\
a_{1,1}\cd a_{1,n}\ovl a_{1,n-1}\cd \ovl a_{1,1}v_{\ovl 1}
+w&C_n,\\
a_{1,1}\cd a_{1,n-1}\ovl a_{1,n-1}\cd \ovl a_{1,1}v_{\ovl 1}
+w&D_n,
\end{cases}
\label{c-L1}
\end{equation}
where $w$ is a linear combination of vectors with higher weights
than the weight of $v_{\ovl 1}$.
Since $\lan v_{\ovl 1},v_{\ovl 1}\ran=1$, 
the coefficient of $v_{\ovl 1}$ in (\ref{c-L1})
is equal to $F^{(n)}_1(y)$ 
and it coincides with the formula in Proposition \ref{Fi} for $i=1$.

In order to show the proposition for the cases $i>1$, 
we need the following lemmas.
\begin{lem}
\label{lem-x}
Let $x_i(c)\in U$ ($i\in I$, $c\in\bbC^\times$) be as above and 
$v_{\Lm_i^{(n)}}$ be the lowest weight vector in 
$L(\Lm_i^{(n)})$ $(i>1)$.
Then we have,
\begin{eqnarray}
&&B_n\text{-case}:\q
x_1(b_1)\cd x_{n-1}(b_{n-1})x_n(a_n)\cd x_1(a_1)v_{\Lm_i^{(n)}}\\
&&\q\qq=\begin{cases}
(b_1\cd b_{i-1})^2(b_i\cd b_{n-1})a_n^2(a_{n-1}\cd a_i)
(+,v_{\Lm_{i-1}^{(n-1)}})+w&\text{ if }i<n,\\
b_1\cd b_{n-1}a_n(+,v_{\Lm_{n-1}^{(n-1)}})+w&\text{ if }i=n,
\end{cases}\nn\\
&&C_n\text{-case}:\q
x_1(b_1)\cd x_{n-1}(b_{n-1})x_n(a_n)\cd x_1(a_1)v_{\Lm_i^{(n)}}\\
&&=(b_1\cd b_{i-1})^2(b_i\cd b_{n-1})(a_na_{n-1}\cd a_i)
\left(\sum_{j=1}^i(-1)^{j-1}v'[v_1;j]\right)+w,
\nn\\
&&
D_n\text{-case}:\q
x_1(b_1)\cd x_{n-2}(b_{n-2})x_n(b_{n-1})x_{n-1}(a_{n-1})
\cd x_1(a_1)v_{\Lm_i^{(n)}}\\
&&\qq=
\begin{cases}
(b_1\cd b_{i-1})^2(b_i\cd b_{n-1})(a_{n-1}\cd a_i)
(+,v_{\Lm_{i-1}^{(n-1)}})+w&\text{ if }i\leq n-2,\\
b_1\cd b_{n-2}a_{n-1}
(+,v_{\Lm_{n-2}^{(n-1)}})+w&\text{ if }i=n-1,\q n:\text{ even},\\
b_1\cd b_{n-2}b_{n-1}
(+,v_{\Lm_{n-2}^{(n-1)}})+w&\text{ if }i=n-1,\q n:\text{ odd},\\
b_1\cd b_{n-2}a_{n-1}
(+,v_{\Lm_{n-1}^{(n-1)}})+w&\text{ if }i=n,\q n:\text{ odd},\\
b_1\cd b_{n-2}b_{n-1}
(+,v_{\Lm_{n-1}^{(n-1)}})+w&\text{ if }i=n,\q n:\text{ even},
\end{cases}\nn
\end{eqnarray}
where $v'$ is as in Lemma \ref{c+} and 
$w$ is a linear combination of weight vectors with 
lower weights than the one of the leading vectors.
\end{lem}
{\sl Proof.}
We can verify this lemma by the induction on the rank $n$.

Let us see the case $B_n$. Since $e_1v_{\Lm_i^{(n)}}=0$ for $i\ne1$, we have
\[
 x_1(b_1)\cd x_{n-1}(b_{n-1})x_n(a_n)\cd x_1(a_1)v_{\Lm_i^{(n)}}
=x_1(b_1)\cd x_{n-1}(b_{n-1})x_n(a_n)\cd x_2(a_2)v_{\Lm_i^{(n)}}.
\]
Suppose that $1<i<n$.
By Lemma \ref{b+-} we have $v_{\Lm_i^{(n)}}=(-,v_{\Lm_{i-1}^{(n-1)}})$
and $x_i(a)$ ($i\ne1$) does not change the top $-$ in
$v_{\Lm_i^{(n)}}$. 
Then, by the hypothesis of the induction we obtain 
\begin{eqnarray*}
&& x_1(b_1)x_2(b_2)\cd x_{n-1}(b_{n-1})x_n(a_n)\cd x_2(a_2)v_{\Lm_i^{(n)}}\\
&&=x_1(b_1)(-,x_2(b_2)\cd x_{n-1}(b_{n-1})x_n(a_n)\cd x_2(a_2)
v_{\Lm_{i-1}^{(n-1)}})\\
&&=(b_2\cd b_{i-1})^2(b_i\cd b_{n-1})a_n^2(a_{n-1}\cd a_i)
x_1(b_1)(-,(+,v_{\Lm_{i-2}^{(n-2)}})))+w\\
&&=(b_1\cd b_{i-1})^2(b_i\cd b_{n-1})a_n^2(a_{n-1}\cd a_i)
(+,(-,v_{\Lm_{i-2}^{(n-2)}})))+w'\\
&&=(b_1\cd b_{i-1})^2(b_i\cd b_{n-1})a_n^2(a_{n-1}\cd a_i)
(+,v_{\Lm_{i-1}^{(n-1)}}))+w'
\end{eqnarray*}
where $w$ and $w'$ are linear combinations of 
lower weight vectors than the leading term.
The case $i=n$ is obtained similarly. 

Next, let us see the case $C_n$. 
As in the previous case, for $i\ne 1$ we have
\[
 x_1(b_1)\cd x_{n-1}(b_{n-1})x_n(a_n)\cd x_1(a_1)v_{\Lm_i^{(n)}}
=x_1(b_1)\cd x_{n-1}(b_{n-1})x_n(a_n)\cd x_2(a_2)v_{\Lm_i^{(n)}}
\]
By Lemma \ref{c+} we get 
\[
v_{\Lm_{i}^{(n)}}=
\sum_{j=1}^{i}(-1)^{i-j}v'[v_{\ovl 1};j].
\]
Since $x_i(a)v_{\ovl1}=v_{\ovl 1}$ for $i\ne 1$, we have
\[
x_i(a) v_{\Lm_{i}^{(n)}}=
\sum_{j=1}^{i}(-1)^{i-j}(x_i(a)v')[v_{\ovl 1};j]\q(i\ne 1), 
\]
and then 
\begin{eqnarray*}
&&x_1(b_1)\cd x_{n-1}(b_{n-1})x_n(a_n)\cd x_2(a_2)v_{\Lm_i^{(n)}}\\
&&=\sum_{j=1}^{i}(-1)^{j-1}x_1(b_1)\{(x_2(b_2)\cd 
x_{n-1}(b_{n-1})x_n(a_n)\cd x_2(a_2)v')[v_{\ovl 1};j]\}\q(i\ne 1).
\end{eqnarray*}
Applying the induction hypothesis to the index set $\{2,3,\cd,n\}$,
we obtain 
\begin{eqnarray*}
&& x_2(b_2)\cd 
x_{n-1}(b_{n-1})x_n(a_n)\cd x_2(a_2)v'\\
&&=(\sum_{j=1}^{i-1}(-1)^{j-1}
(b_2\cd b_{i-1})^2(b_i\cd b_{n-1})(a_na_{n-1}\cd a_i)
v''[v_1;j])+w''
\end{eqnarray*}
where $v''$ is the vector obtained by replacing $v_{\ovl j}$ with 
$v_{\ovl{j+2}}$ in $v_{\Lm_{i-2}^{(n-2)}}$ 
and $w''$ is a linear combination of lower weight
vectors than the ones of the leading term.
Thus, we have
\begin{eqnarray} 
&&\hspace{30pt}x_1(b_1)\cd x_{n-1}(b_{n-1})x_n(a_n)\cd x_2(a_2)v_{\Lm_i^{(n)}}
\label{c1n}\\
&&=(b_2\cd b_{i-1})^2(b_i\cd b_{n-1})(a_na_{n-1}\cd a_i)
\sum_{j=1}^{i}(-1)^{i-j}x_1(b_1)
\left(\sum_{k=1}^{i-1}(-1)^{k-1}
(v''[v_2;k])[v_{\ovl 1},j]\right)+w\nn
\end{eqnarray}
where $w$ is a linear combination of lower weight vectors.
Since $v''$ does not include $v_2$ or $v_{\ovl 1}$, we have
\begin{equation}
 x_1(b_1)(v''[v_2;k])[v_{\ovl 1},j]
=b_1^2(v''[v_1;k])[v_{\ovl 2},j]+\text{lower terms}.
\label{c121}
\end{equation}
Here we have
\[
( v''[v_1;k])[v_{\ovl 2};j]=\begin{cases}
( v''[v_{\ovl 2};j-1])[v_1;k]&\text{ if }k<j,\\
( v''[v_{\ovl 2};j])[v_1;k+1]&\text{ if }k\geq j.
\end{cases}
\]
Thus, 
\begin{eqnarray*}
&&\sum_{j=1}^{i}(-1)^{i-j}
\left(\sum_{k=1}^{i-1}(-1)^{k-1}
(v''[v_1;k])[v_{\ovl 2},j]\right)\\
&&=\sum_{1\leq k<j\leq i}(-1)^{i-j+k-1}( v''[v_{\ovl 2};j-1])[v_1;k]+
\sum_{1\leq j\leq k<i}(-1)^{i-j+k-1}( v''[v_{\ovl 2};j])[v_1;k+1]\\
&&=\sum_{1\leq k\leq j< i}(-1)^{i-j+k}( v''[v_{\ovl 2};j])[v_1;k]+
\sum_{1\leq j<k\leq i}(-1)^{i-j+k}( v''[v_{\ovl 2};j])[v_1;k]\\
&&=\sum_{k=1}^i(-1)^{k-1}\left(\sum_{j=1}^{i-1}(-1)^{i-j-1} 
v''[v_{\ovl 2};j]\right)[v_1;k]=\sum_{k=1}^i(-1)^{k-1}v'[v_1;k].
\end{eqnarray*}
Applying (\ref{c121}) and this to (\ref{c1n}), we obtain
\begin{eqnarray*}
&&x_1(b_1)\cd x_{n-1}(b_{n-1})x_n(a_n)\cd x_2(a_2)v_{\Lm_i^{(n)}}\\
&&=(b_1b_2\cd b_{i-1})^2(b_i\cd b_{n-1})(a_na_{n-1}\cd a_i)
(\sum_{j=1}^i(-1)^{j-1}v'[v_1;j])+w.
\end{eqnarray*}
We have completed the case $C_n$.

Finally, we shall see the case $D_n$.
For $i\leq n-2$ we can show similarly to the case $B_n$. So,
let us show the lemma when $i=n$ and $n$ is odd. In this case,
$v_{\Lm_n^{(n)}}=(-,\cd,-,+)$. 
Since $x_i(a)v_{\Lm_n^{(n)}}=v_{\Lm_n^{(n)}}$ for $i\ne n-1$, 
we have
\begin{eqnarray*}
&& x_1(b_1)\cd x_{n-2}(b_{n-2})x_n(b_{n-1})x_{n-1}(a_{n-1})
\cd x_1(a_1)v_{\Lm_n^{(n)}}\\
&&\qq=x_1(b_1)\cd x_{n-2}(b_{n-2})x_{n-1}(a_{n-1})v_{\Lm_n^{(n)}}.
\end{eqnarray*}
Thus, by direct calculations, we have
\[
x_1(b_1)\cd x_{n-2}(b_{n-2})x_n(b_{n-1})x_{n-1}(a_{n-1})
\cd x_1(a_1)v_{\Lm_n^{(n)}}
=b_1\cd b_{n-2}a_{n-1}(+,v_{\Lm_n^{(n)}})+w.
\]
Other cases are also proved similarly.\qed

Now, we continue the proof of Proposition \ref{Fi}. 
For the case $1<i\leq n$, in order to obtain the 
explicit form of the function $F^{(n)}_i$ we adopt the induction on $n$. 
First, let us see the case $B_n$.
The induction hypothesis for the index set $\{2,\cd,n\}$ and the remark 
in \ref{Bn} mean that 
for $y'=
(y_2(a_{2,2})\cd y_2(\ovl a_{2,2}))\cd
(y_{n-1}(a_{n-1,n-1})y_n(a_{n-1,n})
y_{n-1}(\ovl a_{n-1,n-1}))y_n(a_{n,n})$
\begin{equation}
y'(+,u_{\Lm_{i-1}^{(n-1)}})=
\begin{cases}\left(\displaystyle
\prod_{2\leq k\leq i\leq j\leq n}
a_{k,j}\ovl a_{k,j}\prod_{2\leq k\leq j<i}
\ovl a_{k,j}^2\right)
(+,v_{\Lm_{i-1}^{(n-1)}})+w
&\text{if }i<n,\\
\displaystyle
\prod_{2\leq k\leq j\leq n}
\ovl a_{k,j}(+,v_{\Lm_{n-1}^{(n-1)}})+w&\text{if }i=n,
\end{cases}
\label{bn-1}
\end{equation}
where $\ovl a_{k,n}=a_{k,n}$, 
$w$ is a linear combination of weight vectors with 
higher weight than the one of the leading term. 
We shall denote the
coefficient of $(+,v_{\Lm_{i-1}^{(n-1)}})$ in 
(\ref{bn-1}) by $\Xi_{i-1}^{(n-1)}$.
Set $y:=(y_1(a_{1,1})\cd y_1(\ovl a_{1,1}))y'$. Then we have 
\begin{equation}
\lan yu_{\Lm_{i}^{(n)}},v_{\Lm_{i}^{(n)}}\ran
=\lan y'u_{\Lm_{i}^{(n)}},
x_1(\ovl a_{1,1})\cd x_1(a_{1,1})v_{\Lm_{i}^{(n)}}\ran.
\label{yx}
\end{equation}
By Lemma \ref{lem-x}, (\ref{bn-1}) and the fact that 
$\lan (+,v_{\Lm_{i-1}^{(n-1)}}), (+,v_{\Lm_{i-1}^{(n-1)}})\ran=1$, we have 
\begin{eqnarray*}
&& \lan y'u_{\Lm_{i}^{(n)}},
x_1(\ovl a_{1,1})\cd x_1(a_{1,1})v_{\Lm_{i}^{(n)}}\ran\\
&&\qq =\lan\Xi_{i-1}^{(n-1)}(+,v_{\Lm_{i-1}^{(n-1)}})+w, 
\Omega_{i-1}^{(n-1)}(+,v_{\Lm_{i-1}^{(n-1)}})+w'\ran
=\Xi_{i-1}^{(n-1)}\Omega_{i-1}^{(n-1)}, 
\end{eqnarray*}
where $\Omega_{i-1}^{(n-1)}$ is the coefficient of
$(+,v_{\Lm_{i-1}^{(n-1)}})$ in Lemma \ref{lem-x} and 
$w$ (resp. $w'$) is a linear combination of weight vectors with 
higher (resp. lower) weight than the ones of
$(+,v_{\Lm_{i-1}^{(n-1)}})$.
It is easily to see that 
$\Xi_{i-1}^{(n-1)}\Omega_{i-1}^{(n-1)}$ coincides with $F^{(n)}_i(y)$ 
for type $B_n$ in Proposition \ref{Fi}.

The type $C_n$-case and $D_n$-case with $i\ne n-1,n$ are also done similarly.
Thus, let us see  $F^{(n)}_{n-1}$ and $F^{(n)}_n$ for type $D_n$.
Suppose that $n$ is even.
The induction hypothesis for the index set $\{2,\cd,n\}$ and the remark 
in \ref{Dn} mean that 
we have 
\begin{equation}
y'(+,u_{\Lm_{n-2}^{(n-1)}})= \prod_{2\leq k\leq j<n-1}
\ovl a_{k,j}\cdot(\ovl a_{2,n-1}\cd
a_{n-1,n-1})(+,v_{\Lm_{n-2}^{(n-1)}})+w,
\label{yp}
\end{equation}
for $y'=
(y_2(a_{2,2})\cd y_2(\ovl a_{2,2}))\cd
(y_{n-1}(a_{n-1,n-1})y_{n}(\ovl a_{n-1,n-1}))$, 
where $w$ is a higher term as above and we denote the coefficient 
of $(+,v_{\Lm_{n-2}^{(n-1)}})$ in (\ref{yp}) by $\Xi$.
By Lemma \ref{lem-x}, we have 
\[
 x_1(\ovl a_{1,1})\cd x_1(a_{1,1})v_{\Lm_{n-1}^{(n)}}
=\ovl a_{1,1}\cd\ovl a_{1,n-2}a_{1,n-1}(+,v_{\Lm_{n-2}^{(n-1)}})+w'
\]
where $w'$ is a lower term as above and we denote the coefficient 
of $(+,v_{\Lm_{n-2}^{(n-1)}})$ by $\Omega$. 
Then arguing as above,
we have $F_{n-1}^{(n)}(y)=\Xi\Omega$ for 
$y=y_1(a_{1,1})\cd y_1(\ovl a_{1,1})y'$.
The other cases $F_{n-1}^{(n)}$ ($n$:odd) and $F_n^{(n)}$ are 
showed  similarly.
\qed

\renewcommand{\thesection}{\arabic{section}}
\section{Isomorphisms}
\setcounter{equation}{0}
\renewcommand{\theequation}{\thesection.\arabic{equation}}

Let $\io_0=i_1,\cd,i_L$ be a reduced longest word of $\ge$, 
$B^-_{\io_0}$ as in 
Sect. \ref{sch} and $U^-_{\io_0}$ as in Sect.\ref{u-c}.
\begin{thm}\label{iso}
In case $\ge=A_n,B_n,C_n, D_n$, we have the isomorphism of geometric 
crystals $B^-_{\io_0}\cong U^-_{\io_0}$ by the rational map $\Phi$:
\begin{eqnarray*}
\Phi:&B^-_{\io_0}&\longrightarrow 
\q U^-_{\io_0}\\
&Y_{\io_0}(A_1,\cd,A_L)&\mapsto
\q y_{\io_0}(a_1,\cd,a_L)
\end{eqnarray*}
where $a_j=\Phi_j(A)=(A_1^{{\bf a}_{i_1,i_j}}\cd 
A_{j-1}^{{\bf a}_{i_{j-1},i_j}}A_j)^{-1}$
for $A=(A_1,\cd,A_L)\in (\bbC^\times)^L$.
\end{thm}
{\sl Proof. }
First, we shall see that $\Phi$ is birational.
For the longest element $w_0\in W$, let $L$ be its length.
For $j,k$ with $1\leq j<k\leq L$, set
\[
 P_{j,k}:=\{m=(m_1,\cd,m_t)\in \bbZ^t|j<m_1<\cd<m_t<k,\,\,0\leq t\leq k-j\}.
\]
For $m=(m_1,\cd,m_t)\in P_{j,k}$, set $l(m):=t$.
We understand that $m=\emptyset$ if $t=0$.
For $j,k$ with $1\leq j<k\leq L$, set
\[
 M_{j,k}:=\sum_{m\in M_{j,k}}(-1)^{l(t)}{\bf a}_{i_j,i_{m_1}}
{\bf a}_{i_{m_1},i_{m_2}}\cd {\bf a}_{i_{m_t},i_k},
\]
where ${\bf a}_{i,j}$ is an $(i,j)$-entry of the Cartan matrix.
Let $\cA=({\mathfrak a}_{p,q})_{p,q=1,\cd,L}$ be an integer matrix defined by
\[
 {\mathfrak a}_{p,q}:=\begin{cases}
-1&\text{ if }p=q,\\
-{\bf a}_{i_p,i_q}&\text{ if }p<q,\\
0&\text{ if }p>q.
\end{cases}
\]
It is trivial that the matrix $\cA$ is invertible and its inverse 
is also an integer matrix, denoted by $\cB=({\bf
b}_{p,q})_{p,q=1,\cd,L}$, 
which is indeed given by
\[
 {\bf b}_{j,k}:=\begin{cases}
-1&\text{ if }p=q,\\
M_{i_p,i_q}&\text{ if }p<q,\\
0&\text{ if }p>q.
\end{cases}
\]
Then, we easily know that the following rational map $\Psi$
is an inverse of $\Phi$:
\begin{eqnarray}
\Psi:&U^-_{\io_0}&\longrightarrow 
\q B^-_{\io_0}\label{phi-1} \\
&y_{\io_0}(a_1,\cd,a_L)&\mapsto
\q Y_{\io_0}(A_1,\cd,A_L)\nn
\end{eqnarray}
where $A_j=\Psi_j(a)=a_1^{{\bf b}_{1,j}}\cd 
a_{j-1}^{{\bf b}_{j-1,j}}a_j^{-1}$ for 
$a=(a_1\cd,a_L)\in \bbC^L$.

Next, let us see $\vep_i(\Phi_{\io_0}(Y(A)))=\vep_i(Y_{\io_0}(A))$ for 
$Y_{\io_0}(A)=Y_{\io_0}(A_1,\cd,A_L)$.
As in \ref{expli-UGC}, we have 
\begin{equation}
\vep_i(y_{\io_0}(a_1,\cd,a_L))=
\sum_{1\leq j\leq L,\,i_j=i}a_j.
\end{equation}
We also have the explicit form of $\vep_i(Y_{\io_0}(A))$ as in (\ref{sch-ep}),
\begin{equation}
\vep_i(Y_{\io_0}(A_1,\cd,A_L))=
\sum_{1\leq j\leq L,i_j=i}
\frac{1}{A_1^{{\bf a}_{i_1,i}}\cd A_{j-1}^{{\bf a}_{i_{j-1},i}}A_j}
=\sum_{1\leq j\leq L,i_j=i} \Phi_j(A).
\end{equation}
Thus, by these formula we have 
$\vep_i(\Phi(Y(A_1,\cd,A_L)))=\vep_i(Y(A_1,\cd,A_L))$.

\nd 
Next, let us show that $\Phi\circ e_i^c=e_i^c\circ\Phi$.
Set $e_i^c\circ\Phi(Y_{\io_0}(A))
=y_{\io_0}(a')=y_{\io_0}(a'_1,\cd,a'_L)$. By the formula in 
\ref{subsec-ex-uc}, we obtain
\begin{equation}
a'_j=\frac{\Phi_j(A)}{{L^{(i)}_{m(j)-1}(\Phi(A);c)}^{{\bf a}_{i,i_j}}}
\left(\frac{L^{(i)}_{m(j)-1}(\Phi(A);c)}{L^{(i)}_{m(j)}(\Phi(A);c)}\right)
^{\del_{i,i_j}},
\end{equation}
where $\{j_1,\cd,j_l\}$ is same as in \ref{subsec-ex-uc} and 
$m(j)$ is the number $m$ such that $j_{m-1}<j\leq j_m$.

For $Y_{\io_0}(A)=Y_{\io_0}(A_1,\cd,A_L)$ set $e_i^c(Y_{\io_0}(A)):=
Y_{\io_0}(A'_1,\cd, A'_L)$. Each $A'_j$ is given explicitly by
(\ref{eici}). Denote the numerator of (\ref{eici}) by
$Q_j(c_1,\cd,c_m;c)$. Thus, the denominator is $Q_{j-1}(c_1\cd,c_m;c)$
and 
\[
 A'_j=A_j\frac{Q_j(A_1,\cd,A_L;c)}
{Q_{j-1}(A_1,\cd,A_L;c)}.
\]
Note that $Q_j(A;c)=Q_{j-1}(A;c)$ unless $i_j=i$.
Let us calculate $y_{\io_0}(a''):=\Phi(Y_{\io_0}(A'))$.

\nd
The case $i_j=i$ ($j=j_m$):
\begin{eqnarray*}
&&a''_j=\Phi_j(A')=
({A'}_1^{{\bf a}_{i_1,i_j}}\cd 
{A'}_{j-1}^{{\bf a}_{i_{j-1},i_j}}A'_j)^{-1}\\
&&=(A_1^{{\bf a}_{i_1,i_j}}\cd 
A_{j-1}^{{\bf a}_{i_{j-1},i_j}}A_j)^{-1}
\left(\frac{Q_{j_0}(A;c)}{Q_{j_1}(A;c)}\right)^2\cd
\left(\frac{Q_{j_{m-2}}(A;c)}{Q_{j_{m-1}}(A;c)}\right)^2
\left(\frac{Q_{j_{m-1}}(A;c)}{Q_{j_m}(A;c)}\right)\\
&&=\frac{\Phi_j(A)}{Q_{j_{m-1}}(A;c)Q_{j_m}(A;c)},
\end{eqnarray*}
where $Q_{j_0}(A;c)=1$.

The case $i_j\ne i$ ($j_{m-1}<j<j_m$):
\begin{eqnarray*}
&&a''_j=\Phi_j(A')=
({A'}_1^{{\bf a}_{i_1,i_j}}\cd 
{A'}_{j-1}^{{\bf a}_{i_{j-1},i_j}}A'_j)^{-1}\\
&&=(A_1^{{\bf a}_{i_1,i_j}}\cd 
A_{j-1}^{{\bf a}_{i_{j-1},i_j}}A_j)^{-1}
\left(\frac{Q_{j_0}(A;c)}{Q_{j_1}(A;c)}\right)^{{\bf a}_{i,i_j}}\cd
\left(\frac{Q_{j_{m-2}}(A;c)}{Q_{j_{m-1}}(A;c)}\right)^{{\bf a}_{i,i_j}}
=\frac{\Phi_j(A)}{Q_{j_{m-1}}(A;c)^{{\bf a}_{i,i_j}}}.
\end{eqnarray*}
Here, we can easily see that if $i=i_{j_m}$, 
\[
 Q_{j_m}(A;c)=L_{m}(\Phi(A);c).
\]
This means $a_j'=a_j''$ and then $e_i^c\circ\Phi(Y_{\io_0}(A))
=\Phi\circ e_i^c(Y_{\io_0}(A))$.

Finally, let us show that $\gamma_i\circ \Phi(Y_{\io_0}(A))
=\gamma_i(Y_{\io_0}(A))$.
Since for $y\in U^-$ we have
\[
 \gamma_j(y)=\prod_{i=1}^nF^{(n)}_i(y)^{-{\bf a}_{i,j}},
\]
and the explicit form of $F^{(n)}_i(y)$ as in Proposition \ref{Fi}, 
we use case-by-case method for our purpose.

The $B_n$-case:
Denoting  $\Phi(Y_{\io_0}(A))$ by $y_{\io_0}(a)$ we have 
\begin{eqnarray}
&&a_{i,j}=\begin{cases}
\frac{\displaystyle\prod_{k=1}^i A_{k,j-1}\prod_{k=1}^{i-1} A_{k,j+1}
\prod_{k=1}^{i-1} \ovl A_{k,j-1} \prod_{k=1}^{i-1} \ovl A_{k,j+1}}
{\displaystyle A_{i,j}\prod_{k=1}^{i-1}(A_{k,j}\ovl A_{k,j})^2}&
1\leq i\leq j\leq n-1,\label{aij}\\
\frac{\displaystyle \prod_{k=1}^i A_{k,n-1}\prod_{k=1}^{i-1} \ovl A_{k,n-1}}
{\displaystyle A_{i,n}\prod_{k=1}^{i-1} A_{k,n}^2}&1\leq i\leq j=n,
\end{cases}\\
&&\ovl a_{i,j}=
\frac{\displaystyle\prod_{k=1}^i A_{k,j-1}\prod_{k=1}^{i} A_{k,j+1}
\prod_{k=1}^{i-1} \ovl A_{k,j-1} \prod_{k=1}^{i} \ovl A_{k,j+1}}
{\displaystyle \ovl A_{i,j}\prod_{k=1}^{i}A_{k,j}^2\prod_{k=1}^{i-1}
\ovl A_{k,j}^2}\qq
1\leq i\leq j\leq n-1,\label{ovl-aij}
\end{eqnarray}
where $A_{k,j}=1$ for $k>j$ and $A_{k,n}=\ovl A_{k,n}$.
We shall show the explicit form of  $F^{(n)}_i(\Phi(Y_{\io_0}(A)))$
as follows:
\begin{eqnarray}
&& F^{(n)}_i(\Phi(Y_{\io_0}(A)))=\prod_{m=1}^i\frac{1}{A_{m,i}\ovl
 A_{m,i}}
\q(1\leq i\leq n-1),
\label{F-phi1}\\
&&F^{(n)}_n(\Phi(Y_{\io_0}(A)))=\prod_{m=1}^n \frac{1}{A_{m,n}}.
\label{F-phin}
\end{eqnarray}

Set $q_l(a):=(\ovl a_{l,l}\cd \ovl a_{l,i-1})^2
(\ovl a_{l,i}\cd \ovl a_{l,n-1})a_{1,n}^2(a_{l,n-1}\cd a_{l,i})$
$(1\leq l\leq i)$ and we have
\[
F_i^{(n)}(y_{\io_0}(a))=q_1(a)\cd q_i(a)\q(i<n).
\]
By calculating directly, we get
\begin{equation}
q_l(\Phi(A))=
\frac{A_{l,i}\ovl A_{l,i}}{A_{l,l}^2\ovl A_{l,l}^2}
\prod_{k=1}^{l-1}\frac{A_{k,l-1}^2\ovl A_{k,l-1}^2}
{A_{k,l}^2\ovl A_{k,l}^2},
\end{equation}
and then (\ref{F-phi1}). We also get (\ref{F-phin}) 
by the similar way. The $C_n$-case is obtained similarly. 
So let us see the case $i=n-1,n$ of type $D_n$.
If $n$ is odd, by Proposition \ref{Fi}, we have
\[
F_n^{(n)}(y_{\io_0(a)})=\left(\prod_{1\leq k\leq j<n-1}
\ovl a_{k,j}\right)\cdot(a_{1,n-1}\ovl a_{2,n-1}\cd\ovl a_{n-1,n-1}).
\]
Set $a=\Phi(A)$. 
For $1\leq i\leq j\leq n-2$, $a_{i,j}=\Phi(A)_{i,j}$
and $\ovl a_{i,j}=\ovl{\Phi(A)}_{i,j}$ are same as (\ref{aij}) and 
(\ref{ovl-aij}). For $1\leq i\leq j=n-1$, we have 
\begin{equation}
a_{i,n-1}=
\frac{\displaystyle\prod_{k=1}^{i} A_{k,n-2}
\prod_{k=1}^{i-1} \ovl A_{k,n-2}}
{\displaystyle A_{i,n-1}\prod_{k=1}^{i-1}A_{k,n-1}^2},\qq
a_{i,n-1}=
\frac{\displaystyle\prod_{k=1}^{i} A_{k,n-2}
\prod_{k=1}^{i-1} \ovl A_{k,n-2}}
{\displaystyle \ovl A_{i,n-1}\prod_{k=1}^{i-1}\ovl A_{k,n-1}^2}.
\end{equation}
Then, by calculating directly we have
\[
 F^{(n)}_n(\Phi(Y_{\io_0}(A)))=(\ovl A_{1,n-1}\cd 
\ovl A_{n-1,n-1})^{-1}.
\]
The case for even $n$ is obtained similarly. We have 
\[
 F^{(n)}_{n-1}(\Phi(Y_{\io_0}(A)))=(A_{1,n-1}\cd 
A_{n-1,n-1})^{-1}.
\]
Since 
$\gamma_{n-1}(Y_{\io_0}(A))=(A_{1,n-1}\cd 
A_{n-1,n-1})^{-1}$ and 
$\gamma_{n}(Y_{\io_0}(A))=(\ovl A_{1,n-1}\cd 
\ovl A_{n-1,n-1})^{-1}$, 
we know that $\gamma_i(Y_{\io_0}(A))
=\gamma_i(\Phi(Y_{\io_0}(A)))$ and completed 
the proof of Theorem \ref{iso}.\qed


Since $B^-_{\io_0}$ (resp. $U_{\io_0}^-$) is birationally 
equivalent to the flag variety $X$ 
(resp. unipotent radical $U^-\subset B^-$), we have the following:
\begin{cor}
There exists an isomorphism of geometric crystals:
$X\cong U^-$.
\end{cor}

\renewcommand{\thesection}{\arabic{section}}
\section{Conjectures}
\setcounter{equation}{0}
\renewcommand{\theequation}{\thesection.\arabic{equation}}

For a Weyl group element $w\in W$, let $\io={i_1}\cd {i_k}$
be a reduced word of $w$.
Set 
$U_{\io}^-:=\{y_{i_1}(c_1)\cd y_{i_k}(c_k)|
c_1,\cd,c_k \in\bbC^\times\}$. 
Let  $u_w^{(i)}$ be the normalized extremal weight vector
with the extremal weight $w\Lm_i$
in $L(\Lm_i)$ and 
define the function on $U_{\io}^-$ by 
\[
 F^{(n)}_{(i,\io)}(y):=\lan y\cdot u_{\Lm_i},u_w^{(i)}\ran\q
(y\in U_{\io}^-,\,\,i\in I).
\]
Here
we present the following conjecture:\\
\vskip0mm
\nd{\bf Conjecture.}
{\it If $G$ is semi-simple and $I=I(w)$ 
(see Sect.\ref{sch}.), we can associate a
 geometric crystal structure 
with $U_{\io}^-$ and it is isomorphic to the geometric crystal on 
the Schubert variety $\ovl X_w$ for any $w\in W$.}

For a reduced word $\io$, set 
$B_{\io}^-:=\{Y_{i_1}(c_1)\cd Y_{i_k}(c_k)|
c_1,\cd,c_k \in\bbC^\times\}$ 
 ($Y_i(c)=y_i(\frac{1}{c})\al_i^\vee(c)$) and
let 
$(B_{\io}^-,\{e_i\}, \{\gamma_i\},\{\vep_i\})$
be the geometric crystal isomorphic to $\ovl X_w$
as in Sect.\ref{sch} (see also \cite{N}). 
To show the conjecture, we should obtain that 
$U_{\io}^-\cong B_{\io}^-$ as geometric crystals.

\end{document}